\documentclass[preprint,10pt]{elsarticle}

\usepackage{amssymb}
\usepackage{amsmath}
\usepackage{lineno}
\usepackage{algorithm}
\usepackage{algorithmic}
\usepackage{subcaption}
\usepackage{relsize}
\usepackage[utf8x]{inputenc}
\graphicspath{ {./pics/} }
\usepackage[title]{appendix}
\journal{arXiv.org}
\newcommand{\pd}  {\partial}
\newcommand{\dext}{\mathrm{d}}
\renewcommand{\vec}[1]{\mathbf{#1}}
\newcommand{\norm}[1]{\left\lVert#1\right\rVert}
\renewcommand{\i}{{\mathrm{i}}}

\begin{document}
\begin{frontmatter}

\title{Quad layouts with high valence singularities\\for flexible quad meshing}
\author{Jovana Jezdimirovi\' c, Alexandre Chemin, Maxence Reberol,\\
	Fran\c{c}ois Henrotte, Jean Fran\c{c}ois Remacle}
\address{Universit\' e catholique de Louvain, Louvain la Neuve, Belgium jovana.jezdimirovic@uclouvain.be}

\begin{abstract}
A novel algorithm that produces a quad layout based on imposed
set of singularities is proposed. 
In this paper, we either use singularities that
appear naturally, e.g., by minimizing Ginzburg-Landau energy, or
use as an input user-defined singularity pattern,
possibly with high valence singularities that do not appear naturally
in cross-field computations. 
The first contribution of the paper is the development of
a formulation that allows computing a cross-field from a given set of
singularities through the resolution of two linear PDEs. A specific
mesh refinement is applied at the vicinity of singularities 
to accommodate the large gradients of cross directions that appear
in the vicinity of singularities of high valence.
The second contribution of the paper is a correction scheme that repairs limit cycles and/or non-quadrilateral patches. Finally, a high quality block-structured quad mesh is generated from the quad layout and per-partition parameterization.
\end{abstract}


\end{frontmatter}

\section{Introduction}
\label{S:Intro}
Quad meshing is a discipline that has been developed for different
purposes by two distinct communities. Computer graphics and
engineering analysis communities have indeed produced
extensive research in the subject for decades, leading to 
both unstructured and structured approaches of quad meshing.
A recent survey \cite{Campen:2017} comprehensively
discusses the issue of quad meshing from manual to fully automatic
generation, along with the respective advantages and drawbacks of the
different methods.

\begin{figure*} 
\begin{center}
  \includegraphics[width=0.95\linewidth]{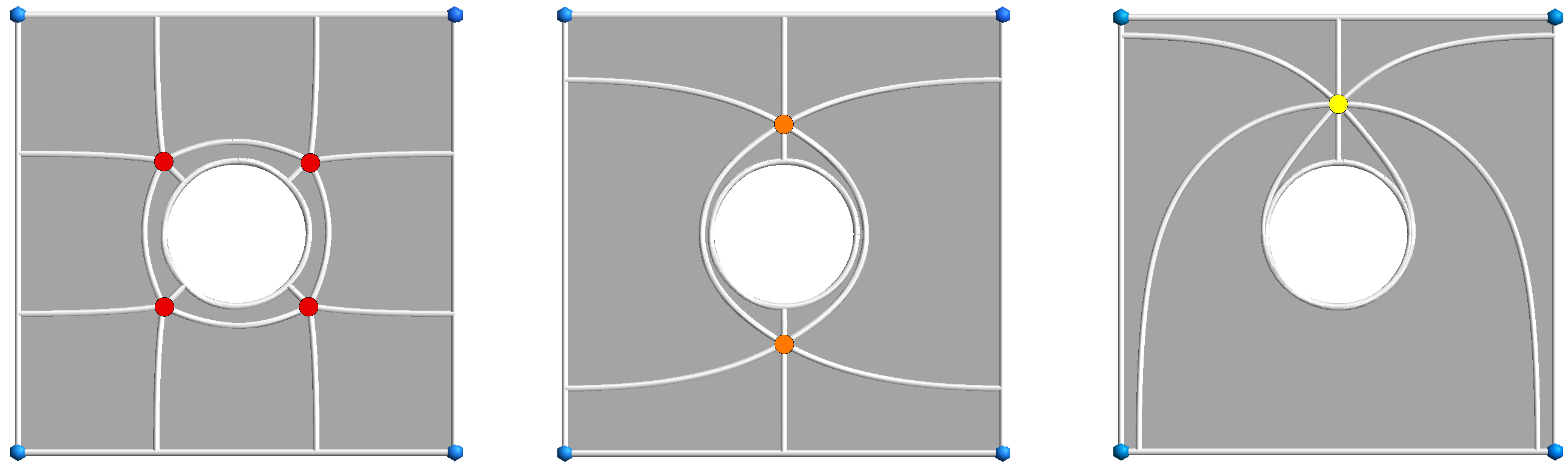}
  \caption{Three quad layouts of a simple domain defined by a square
    minus a disk. Colored points are the irregular corners of the layouts.}
  \label{fig:teaser}
\end{center}
\end{figure*}

In a block-structured quad mesh, most of
the vertices are regular. A vertex $v$ is regular if it has the optimal
number $n_v$ of adjacent quads: $n_v=2$ if the vertex is on the
boundary $\partial \Omega$ of
the domain $\Omega$ to be meshed, and $n_v=4$ if it is an internal vertex.
Our approach for block-structured quad meshing relies on the computation of an
auxiliary object, called a cross-field, as a guide to
orient the quadrangular elements. 
Cross-fields have isolated singularities that
can be shown to correspond to the location of irregular vertices of the quad mesh \cite{Beaufort:2017}. 
Using this useful structure
of cross-fields has become popular to generate block-structured quad
meshes \cite{Bommes:2012, Kowalski:2013, Fogg:2015, Viertel:2019,  Jezdimirovic:2019}.
More specifically, generating block-structured quad meshes requires the
construction of a \emph{quad layout}.  A quad layout is the
partitioning of an object's surface into simple networks of conforming quadrilateral
patches \cite{Campen:2014}. Fig.~\ref{fig:teaser} shows examples of quad
layouts for a simple 2D domain.

Two main approaches exist to extract a quad layout from the
singular structure of a cross-field.
The first approach consists in computing a seamless global
parameterization of the domain where integer iso-values of the
parameter fields form the sides, or the separatrices,
of the quad layout  \cite{Campen:2015, Bommes:2013, Ray:2006}. This approach is implicit in the sense that it does not explicitly compute separatrices from the cross-field. The second approach is explicit. It is the one used in this paper. It consists in directly tracing the separatrices of the quad layout, starting at the singularities of the cross-field. A singularity of valence $5$ corresponds to a vertex $v$ with $n_v=5$ adjacent quads,
so that $5$ separatrices can be started at this vertex $v$. 
Separatrices then follow integral lines of the cross-field, eventually ending on $\partial \Omega$ or on another singularity.

The main drawback of the explicit approach is related to
its lack of geometrical robustness. 
Assuming that
singularities are located at vertices of our mesh,
the orientation of the cross-field is thus varying abruptly at the vicinity
of a singular vertex $v$. That bad representation of the cross-field
can have dramatic consequences.
In many cases, the right number of separatrices cannot even be traced starting at $v$.
Bad cross representation around singularities is also the cause of
tangential crossings, i.e., separatrices that should meet at a singular
vertex but that miss the rendez-vous and continue further, 
eventually forming an infinite loop.

\begin{figure*}[tb]
\begin{center}
\includegraphics[width=0.6\linewidth]{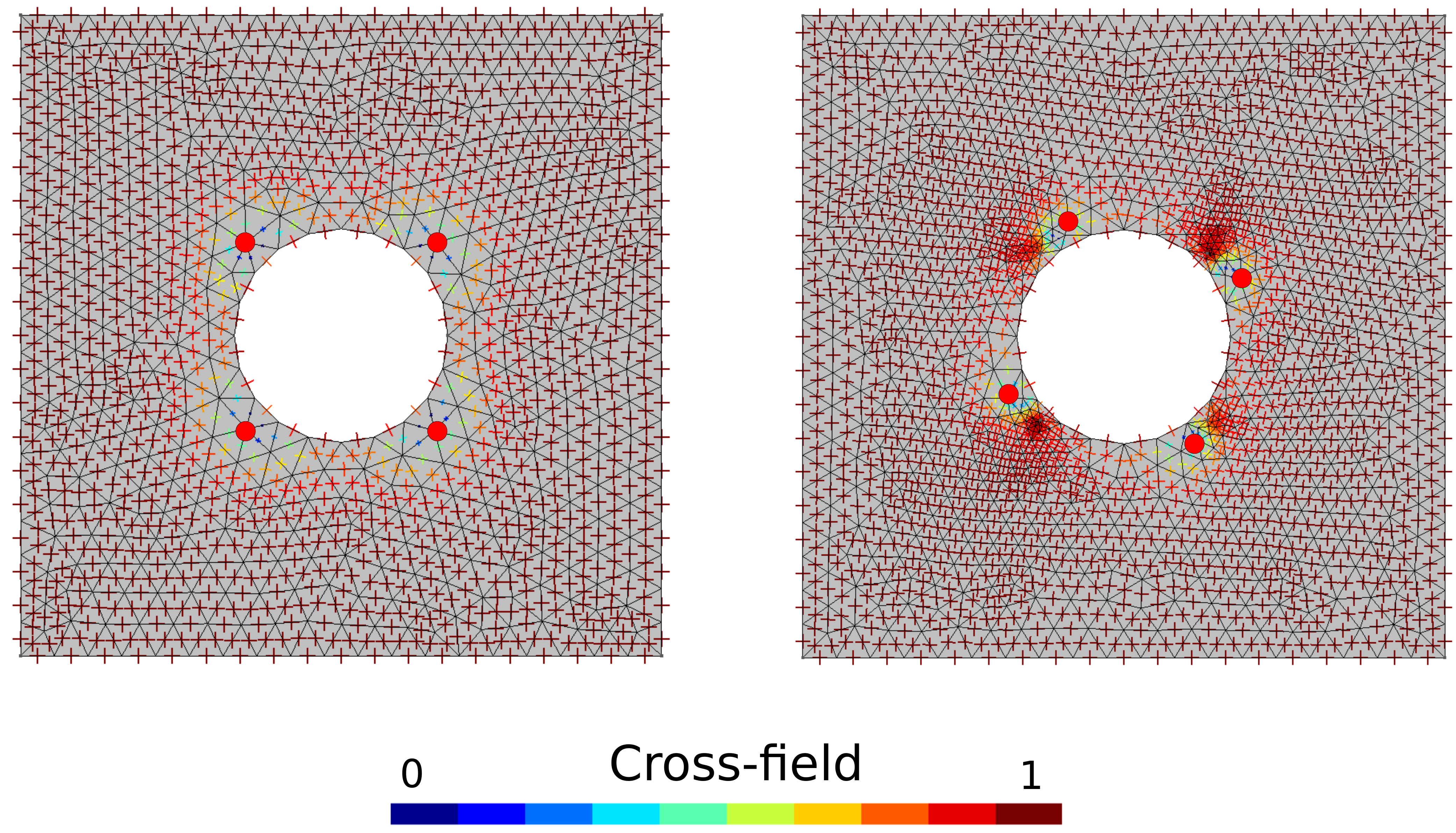} 
\caption{Displacement of singular points due to a cross-field mesh dependency. 
Cross-field computation follows the work from \cite{Beaufort:2017}.}
\label{fig:GL-GL}
\end{center}
\end{figure*}

Mesh refinement around singularities would then seem to be a
reasonable idea. Disappointingly, it is not.
Cross-field smoothers have the bad habit to push singularities away from refined zones, 
because the energies used in cross-field models are unbounded 
(in the continuous case) at singularities,
so that locating singularities in coarse region allows finding better 
minima (in the discrete case). 
Standard mesh adaptation thus fails here, because the
solution does not converge in terms of energy. 
Fig.~\ref{fig:GL-GL} shows an example of that typical behavior. 
Fig.~\ref{fig:GL-GL} (left) shows a cross-field computed on a uniform mesh. 
Mesh refinement is then applied close to those singularities, 
and the cross-field solver is then applied once again with the refined mesh. 
On Fig.~\ref{fig:GL-GL} (right), one can observe that the singularities have moved with respect to their position after the first cross-field computation. 
The tendency to escape refined regions is so strong that the solution is
actually breaking the symmetry of the problem. 

In order to avoid this issue,
our cross-field computation is implemented as a two-step process.
A first cross-field is computed in a standard fashion by means of
minimizing a nonlinear energy 
(the Ginzburg-Landau energy functional in our case).
The {\em singularity pattern} of the cross-field is extracted from
this first computation and then injected as a constraint for a second
cross-field computation on an adapted mesh. The formulation for
computing this second cross-field is different. Here, the singularity
pattern is taken as an input and a cross-field is computed by means of
solving two linear systems. Singularity positions \emph{are thus prescribed} so mesh adaptation can be performed.
The local structure of a cross-field at a singularity is essentially
radial: the mesh that is used in this second step
thus involves specific \emph{bicycle spokes patterns} at singularities (see
section~\ref{S:singularities}). With that {\it ad hoc}
refinement, the singular structure of the cross-field is perfectly
captured. Another advantage of our approach is the possibility of
moving singularities, adding or removing some, while still respecting
topological constraints.
In this work, we use that advantage for correcting or adapting the
quad layout when necessary, e.g., to fix non-quadrilateral partitions
or to avoid limit cycles. 

\begin{figure*}[h!t]
\begin{center}
\includegraphics[width=\linewidth]{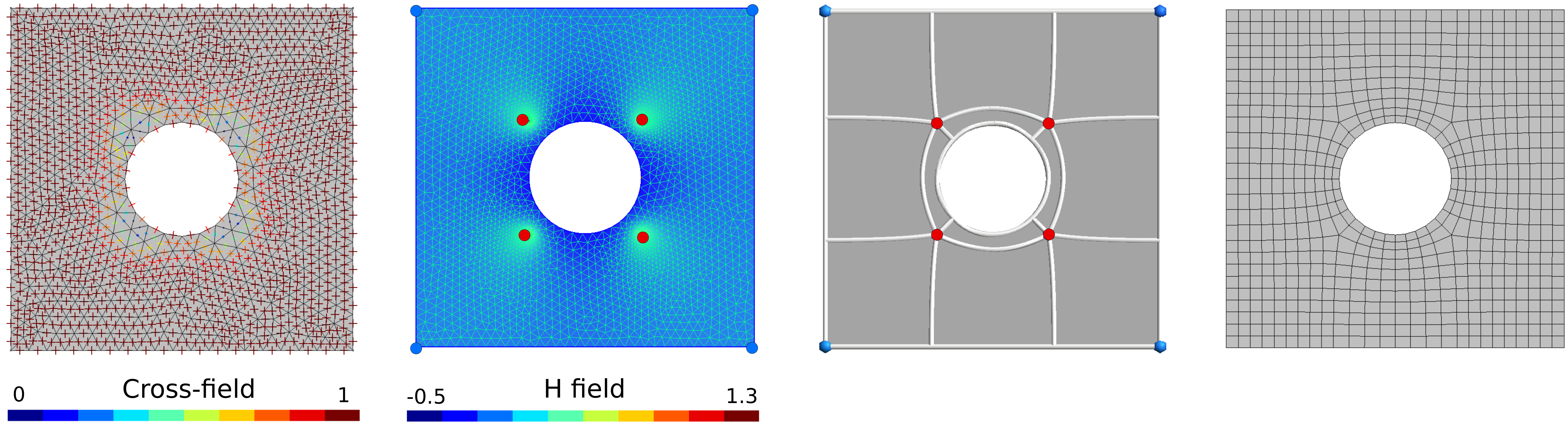} 
\caption{Illustration of 4-steps pipeline}
\label{fig:pipeline}
\end{center}
\end{figure*}

The full pipeline of our block-structured quad mesher consists in four steps
and is described in Figure~\ref{fig:pipeline}.

\noindent{\bf Step 1}: compute a first cross-field using standard non
linear minimization procedure that allows to find a singularity
pattern i.e. singularity positions and indices.

\noindent{\bf Step 2}: compute a second cross-field with the prescribed
  singularity pattern of Step 1 on a refined mesh.

\noindent{\bf Step 3:} compute a neat quad layout decomposition
without infinite loops or t-junctions. Quad layout partitioning is
performed on the accurate cross-field of Step 2
by applying the separatrice tracing scheme presented in \cite{Jezdimirovic:2019}. 

\noindent{\bf Step 4:} generate a full quadrilateral mesh.

\begin{figure*}[h!t]
\begin{center}
\includegraphics[width=0.85\linewidth]{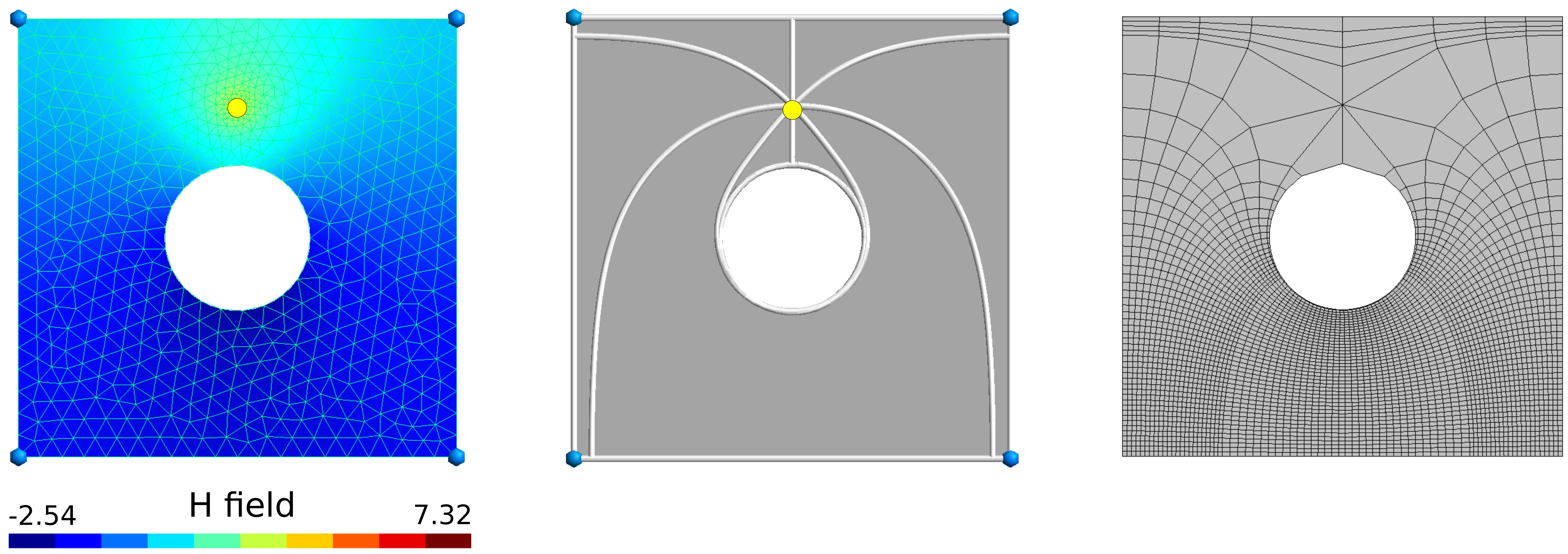} 
\caption{Illustration of the reduced 3-steps pipeline}
\label{fig:pipelineArtifical}
\end{center}
\end{figure*}

Note that Step 1 can be skipped if a singularity pattern is provided by
the user.  Figure \ref{fig:pipelineArtifical} shows this reduced 3-step
process with a very unique singularity of very high valence $8$!

The two-step cross-field computation delivers thus an accurate mesh independent cross-field, Fig.~\ref{fig:pipeline}, from which one can then continue with a neat quad layout decomposition, and eventually a full quadrilateral mesh, Fig.~\ref{fig:pipeline} and \ref{fig:pipelineArtifical}. 
Quad layout partitioning is performed on the final (accurate) cross-field by applying the separatrice tracing scheme presented in \cite{Jezdimirovic:2019}. 
The final step is the generation of the quadrilateral mesh by following the parameterization imposed by the cross-field (section~\ref{S:QuadMesh}).

The steps of the proposed approach are detailed in Algorithm~\ref{alg:pipeline1}. 
The whole pipeline has been implemented as a fully functional module in Gmsh \cite{Geuzaine:2009}, the open source finite element mesh generator.
Although the pipeline can be fully automated, a certain level of interactivity is possible. The user is indeed given the opportunity to modify the singularity pattern and to adapt the location and index of singularities, as long as it corresponds to the topology of the domain. The modified singularity valence can be as high as six or eight, Fig.~\ref{fig:pipelineArtifical}.

\begin{algorithm}[H]
\caption{Generate all quad mesh}
  \label{alg:pipeline1}
\vspace*{0.2cm}
\textbf{INPUT:}  FILE in \textbf{.geo} format \\
\textbf{OUTPUT:}  FILE in \textbf{.msh} format 

\begin{algorithmic}[1] 
\STATE Non-linear cross-field computation on $\Omega$ (Sect.~\ref{cf1})
\STATE Extract singularity pattern (Sect.~\ref{S:singularities})
\STATE Compute cross-field with imposed singularities (Sect.~\ref{S:Hfunc})
\STATE Obtain and correct quad layout (Sect.~\ref{S:MultiBlock})
\STATE Generate quad mesh on $\Omega$ (Sect.~\ref{S:QuadMesh})
\end{algorithmic}
\end{algorithm}

\section{Related work}
\label{S:Related}

Due to the fact that some of the techniques used in our approach are customary in the meshing and 
computational geometry communities, we will focus in this article on the most notable trends and contrasts/resemblances with existing works. For the sake of comprehensiveness, bibliographic references for our algorithm are indicated in corresponding sections.

\subsection{Cross-field generation for meshing purposes}

For quad meshing purposes, authors are generally relying on either explicit quadrangulations or parameterization/remeshing techniques \cite{Bommes:2009}. The group of parameterization approaches is enriched with numerous methods \cite{Campen:2017}, where rough dissemination could be made with two generalized flows:
cross-field guidance algorithms (\cite{Kowalski:2013, Fogg:2015, Viertel:2019, Jezdimirovic:2019})
and cut-graph techniques (\cite{Campen:2015, Bommes:2013}). Our approach follows the first flow with the difference of using two independent 
successive cross-fields, which is unusual in related work on cross-field generation (\cite{Ray:2008, Vaxman:2016}). Similarly  to \cite{Jezdimirovic:2019, Ray:2006, Diamanti:2015}, we aim at obtaining a locally integrable cross-fields, such that for each point $X$ inside a patch holds:
\begin{equation}
\nabla \times {\bf u}_1=0, \quad \nabla \times {\bf u}_2=0,
\end{equation}
where ${\bf u}_1$ and ${\bf u}_2$ represent two orthogonal branches of the cross defined in point $\bf X$.

Stressed out by \cite{Crane:2010}, a major compromise to be done in the computation of cross-fields
is between efficient design of vector field and control over singularities. Our algorithm runs fully automatic which frees the user of painful task to explicitly specify all singularities on complex domains \cite{Ray:2009}. At the same time, due to the independent steps of the algorithm, we are also able to fulfill the aim of generating a fast smooth cross-field that respects a user-defined singularity patterns as in \cite{Bommes:2009, Ray:2008}.

\subsection{Singularities with high indices}

High-order singularities' patterns play an important role in numerous scientific applications (\cite{Scheuermann:1997, Scheuermann:1998, Tricoche:2000, Gunther:2018, Weinkauf:2005}). Initially, high-order singularities' definition, extraction and visualization were relying on mathematical background of Clifford algebra, typically  motivated by applications for fluid mechanics or electrostatics (\cite{Scheuermann:1997, Gunther:2018, Weinkauf:2005}). Some more recent papers (\cite{Tricoche:2000, Weinkauf:2005}) are focused on their role in vector field generation. In our case, the application of high-order singular points is exploited for its benefits for quad meshing purposes, specifically obtaining a quad layout and inducing important element size gradients in the refined block-structured quad mesh.

\subsection{Quad layout partition guided by cross-field}

Representation of geometrically or topologically complex domains requires their decomposition into blocks, preferably quadrilateral. A wide range of existing methods is described in the survey \cite{Campen:2017}. 
Works analogous to ours, i.e. cross-field based approaches, are decomposing the domain either by generating separatrice graph obtained from numerical integration of streamlines (\cite{Kowalski:2013, Fogg:2015, Viertel:2019, Jezdimirovic:2019, Viertel2:2019}) or by constructing edge maps (\cite{Myles:2014, Ray:2014, Jadhav:2011}) and  motorcycle graph (\cite{Eppstein:2008}) using the isolines of an underlying parameterization. As already reported by authors working on the first group of approaches, generated quad layout commonly requires repairing of invalid configurations caused by the appearance of limit cycles and thin blocks. Some of the efficient methods for solving these issues are \cite{Myles:2014, Viertel2:2019, Borden:2002, Daniels:2008} as well as the quantization method \cite{Campen:2015}. Generating and correcting the quad layout, in our case, shares some common points with \cite{Jezdimirovic:2019, Campen:2015,  Viertel2:2019}, exact differences with respect to these works will be outlined in the corresponding sections.

\section{First cross-field computation}
\label{cf1}

\begin{figure}[h!t]
  \begin{center}
    \includegraphics[width=0.2\textwidth]{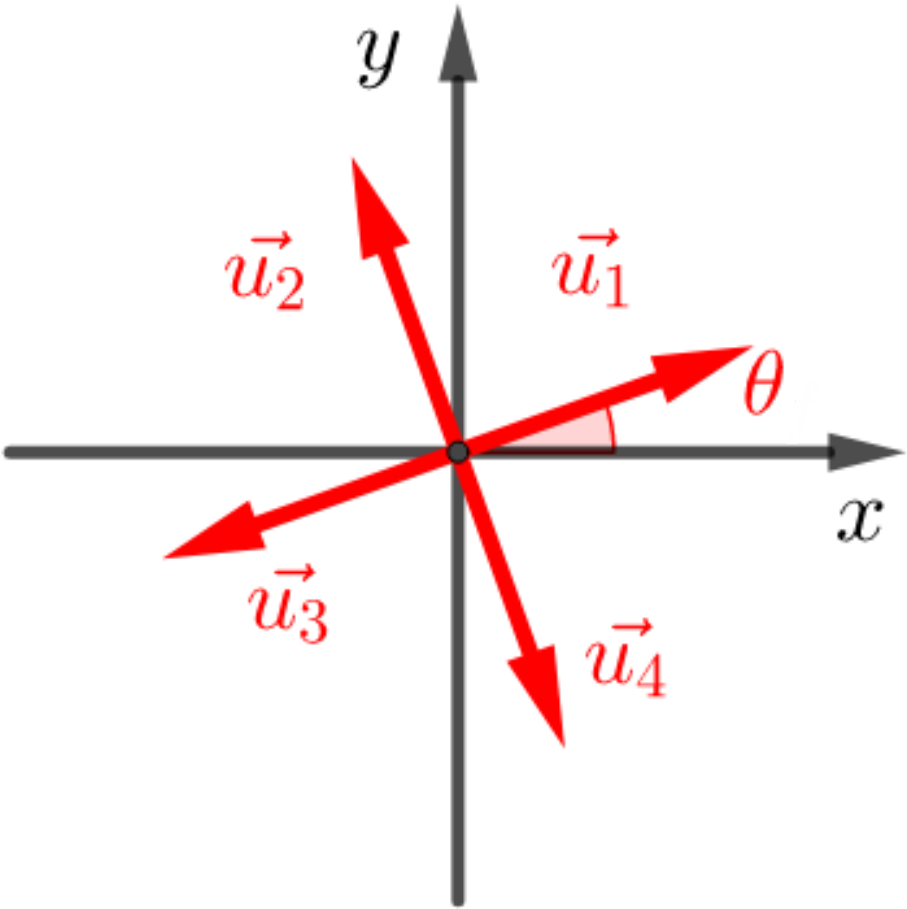}
    \caption{2D cross definition as a set of four orthogonal vectors}
    \label{fig:cross}
  \end{center}
\end{figure}

A 2D cross $\mathbf{c}$ is a set of $4$ orthogonal vectors of the same norm $l$, 
i.e., $\mathbf{c}=\{ \mathbf{u}_k , |\mathbf{u}_k|=1 , k\in [|1,4|] \}$. 
On a planar 2D domain $\Omega$,
a Cartesian basis $\{ \mathbf{x},\mathbf{y} \}$ can be defined,
and the angular orientation of the cross $\mathbf c$ can be simply represented by an angle $\theta$ in that basis,
since one has 
$\mathbf{u}_k=l\cos (\theta + k\frac{\pi}{2})\mathbf{x} + l\sin (\theta + k\frac{\pi}{2})\mathbf{y}$, 
as depicted in Fig.~\ref{fig:cross}.
This angular representation $\theta$ of the cross is however determined 
only up to an additive factor $k' \frac\pi 2$ with $k'\in \mathbb Z$. 

A 2D cross-field, now, is a map $\mathcal{C}_\Omega: \mathbf{X} \in \Omega \to \mathbf{c}(\mathbf{X})$,
and the standard approach to compute a smooth boundary-aligned cross-field is to minimize the Dirichlet energy
\begin{equation}
    \min_{\mathcal{C}_\Omega} \int_\Omega \norm{\nabla \mathcal{C}_\Omega}^2 
    \label{eq:dirichlet_energy}
\end{equation}
subject to the boundary condition $\mathbf{c(\mathbf{X})} = \mathbf{g}(\mathbf{X})$ on $\partial \Omega$,
where $\mathbf{g}$ is a given function. 

The indetermination on $\theta$ mentioned above 
prevents from using it directly as unknown for the finite element formulation of the problem. 
In order to have a unique (symmetry invariant) and continuous cross representation, 
the alternative vector representation 
$\mathbf{v}(\mathbf{X}) = (\cos(4 \theta), \sin(4 \theta))$ is preferred \cite{Hertzmann:2000},
which leaves us with a vector field to solve for.

This vector field is continuous, but also subject to the nonlinear constraint $|\mathbf{v}(\mathbf{X})|=1$.
An elegant way to deal with this constraint is to minimize a Ginzburg-Landau type functional, 
which penalizes local vector values that deviate from the unit norm \cite{Beaufort:2017}. 
Direct solving of the non-linear Ginzburg-Landau problem is however rather expensive.
Another approach is the MBO algorithm proposed in \cite{Viertel:2019}.
It is a faster and simpler alternative that provably
converges to the same results on uniform meshes. 
The main idea of the MBO algorithm is to alternate a heat diffusion problem, 
that performs the Laplacian smoothing of the vector field, 
with a projection back onto the space of unit norm cross-fields, 
until convergence is achieved. 
Each step of MBO algorithm is thus given by
\begin{equation}
    {\frac{\partial \mathbf{v}} {\partial t}} = \alpha \nabla^2 \mathbf{v}
    \; \text{ then } \;
    \mathbf{v} \leftarrow \frac{\mathbf{v}}{\| \mathbf{v} \|}. 
\end{equation}

As in \cite{Palmer:2020}, we start with large diffusion
(e.g. $\alpha = (0.1 * \text{bbox diagonal})^2$) and
we progressively decrease until we reach the small diffusivity coefficient
associated to the mesh size (e.g. $\alpha = (\text{edge min})^2$).
This strategy considerably accelerates
convergence as the low frequencies are solved before the larger ones,
as in a multi-grid approach. To further accelerate the computation, we
use a small number of diffusivity coefficients (e.g., five to ten) and we apply
the MBO step until convergence for each diffusivity coefficient. As the diffusion
linear system is not changing for a given diffusivity coefficient, we can
factorize the linear system one time and reuse it until convergence.
In the end, the computational cost is approximately one linear system
solve per diffusivity coefficient, i.e., five to ten in total, which is very reasonable
to accurately solve the cross-field non-linear problem (Eq.~\ref{eq:dirichlet_energy}).

\section{Singularity pattern}
\label{S:singularities}

Singularities play an important role for meshing purposes and their properties are well documented in the literature \cite{Kowalski:2013, Fogg:2015, Viertel:2019, Jezdimirovic:2019, Viertel2:2019}. 
Obviously, cross-fields have a close parenthood with vector fields
and, like any direction field (a normalized vector field, a cross-field, \dots, more details in \cite{Ray:2008}),
they are subject to the topological constraints of the Poincar\'e-Hopf theorem.
Let us consider, for instance,  a smooth unit norm vector field $\vec u$ on a surface $\Omega$. 
The local value of this vector field 
can be defined as the rotation by a certain angle $\theta$ with respect to a local angular reference 
around the local normal vector to the surface. 
This natural geometric representation of the vector field has however a number of issues.
It is only defined up to an additive angle $2k\pi, k\in \mathbb Z$.
It is also discontinuous wherever the value of $\theta$  jumps from 0 to $2\pi$ or reversely.
Finally, as customary in differential geometry,
$\theta$ values at distant points are incommensurable unless a parallel transport rule is added,
because they refer to different local normal vectors and different local reference positions.

Still, in an infinitesimal neighborhood $V_X$ of a point $X\in \Omega$,
the parallel transport is superfluous because both the normal vector
and the reference angular reference can be considered constant.  
The smooth vector field $\vec u$ then 
induces a map $C^1 \rightarrow C^1$, where $C^1$ is the unit circle, 
for any (infinitesimal) closed curve $C \subset V_X$. 
This map  follows the turns done by $\vec u$ 
as the curve $C$ is travelled once.
As $\vec u$ is smooth, the number of turns is necessarily an integer number,
called the {\em index} of $\vec u$ at $X$, $\mathrm{index}_X(\vec u) \in \mathbb Z$.
The point $X$ is regular if the index is zero,
and it is a singularity of index $\mathrm{index}_X(\vec u)$ otherwise. 

The principle is similar for a cross-field $\mathcal{C}_\Omega$
except that the map is $C^1 \rightarrow C^1/O$
where $O$ stands for the symmetry group of the cross,
and that the index is now an integer multiple of $\frac 1 4$.
If one notes $S_j, j=1, \dots N$, the singularities of the cross-field, 
the Poincar\'e-Hopf theorem states that 
\begin{equation}
\sum_{j=1}^N \mathrm{index}_{\mathcal{C}_\Omega}(S_j) = \chi(\Omega),
\label{eq:Poincare-Hopf}
\end{equation}
where $\chi(\Omega)$ is the Euler characteristic of the (non necessarily planar) surface $\Omega$.
The indices of the singularities will be noted
\begin{equation}
k_j\equiv \mathrm{index}_{\mathcal{C}_\Omega}(S_j) = \frac t 4
\quad , \quad t \in \mathbb Z
\end{equation}
and we shall also use the related concept of valence defined by
$$
k_j = \frac{4-\mathrm{valence}(S_j)}{4}
\quad , \quad 
\mathrm{valence}(S_j) = 4 - 4 k_j.
$$

\begin{figure}[h!t]
\begin{center}
\includegraphics[width=0.55\textwidth]{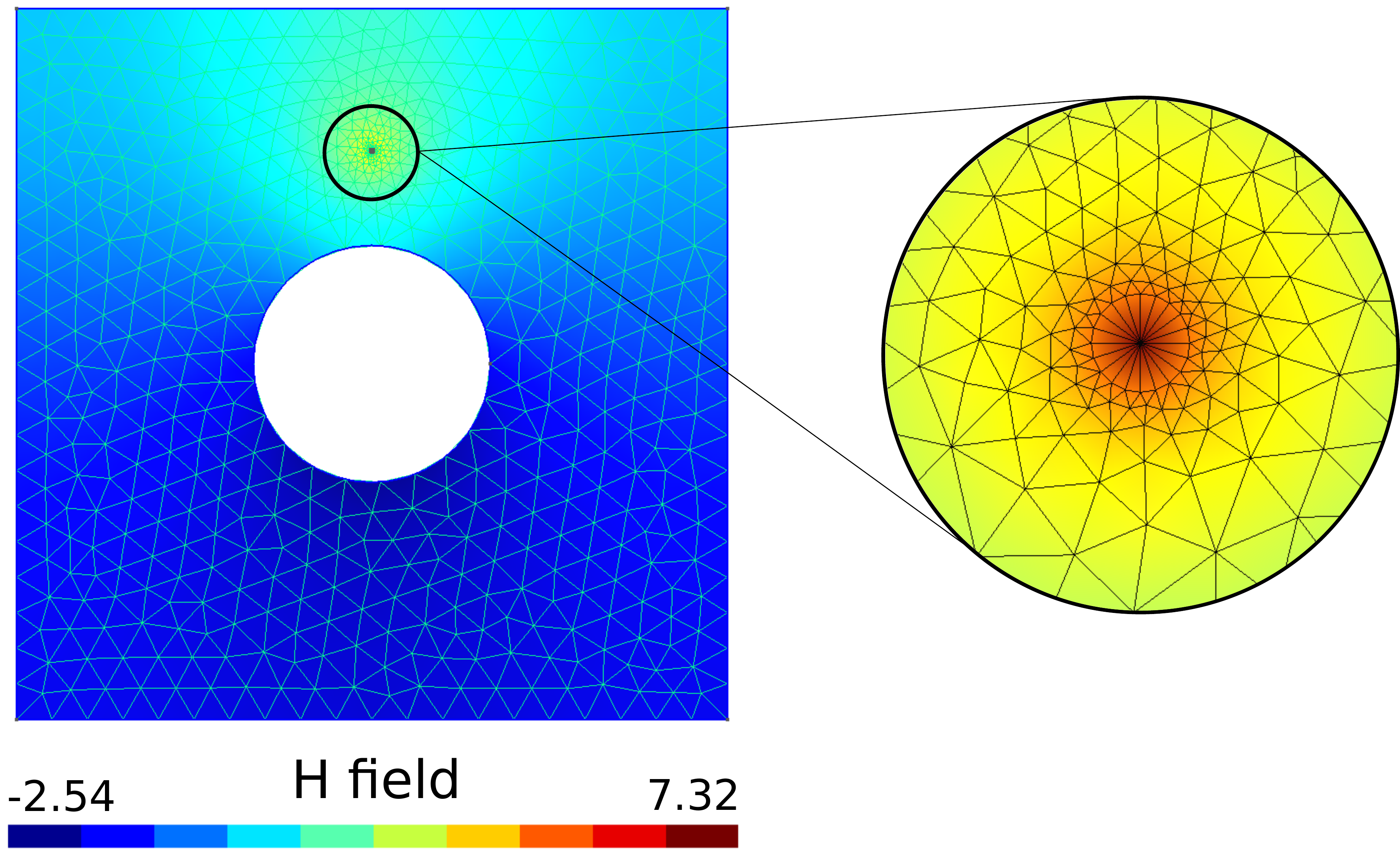} 
\caption{Mesh refinement around singularity of valence 8}
\label{fig:refinement}
\end{center}
\end{figure}

In practice, we use a simple way to find the singularity location  by computing the winding number in a similar fashion as in \cite{Viertel2:2019}, originally from \cite{Kronecker:1869}. The corresponding valence is then determined by detecting the number of separatrices concurring at the singularity \cite{Jezdimirovic:2019}.
When the cross-field is computed on a domain with a coarse mesh around singularities (of valence five or higher),  numerical inaccuracies can be expected in the determination of the cross-field orientations, which will finally lead to an invalid quad layout generation. 
To prevent this, the mesh is refined in the form of bicycle spokes in the vicinity of the identified singularities, so that the large gradients of cross directions are well accommodated, as shown in Fig~\ref{fig:refinement}.

\begin{figure}[tb]
\begin{center}
\includegraphics[width=0.55\textwidth]{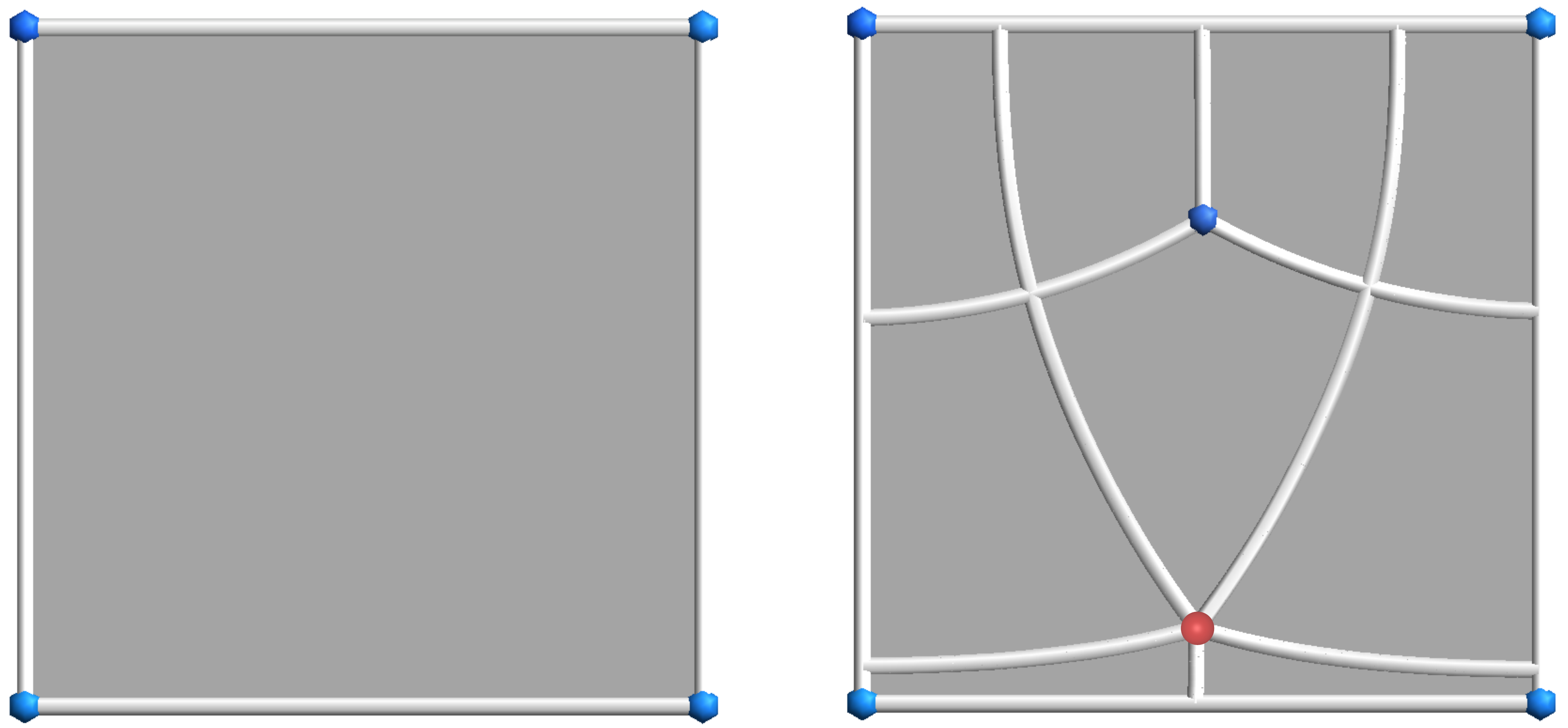} 
\caption{Two quad layouts of the same region. The addition of a 3-5 pair to the singularity pattern,
in the second example, allows imposing a specific size map on the region.}
\label{fig:singType}
\end{center}
\end{figure}

An adaptation of the singularity pattern can be performed at this stage of the pipeline. With the two-step approach, the user is indeed offered the opportunity before the second cross-field computation to modify the location of singularities, or even to place additional singularities, provided the final singularity pattern still fulfills the Poincar\'e-Hopf topological constraint (\ref{eq:Poincare-Hopf}). For instance, additional pairs of singularities of valence $3$ and $5$ do not modify the result of the summation in  (\ref{eq:Poincare-Hopf}). Such 3-5 pairs can therefore be added freely to the singularity pattern in order, for instance, to allow a less distorted quad layout or to impose a specific size map, as shown in Fig.~\ref{fig:singType}.

In practice, imposing a singularity pattern simply consists in encoding valences data for vertices of the initial mesh and  ensuring a sufficient mesh refinement in their respective vicinities.

It will also be shown at the end of the following section that not all valid
singularity patterns will be suitable for obtaining a high quality quad layout,
depending on the geometrical characteristics of the domain.  Mispositioning of
singularities or imposing inadequate valences may lead in practice to
non-quadrilateral patches. 

\section{Cross-field computation with imposed singularity pattern}
\label{S:Hfunc}

The nonlinear problem described in Sect.~\ref{cf1}
is able to reveal a singularity pattern 
for the domain $\Omega$ under analysis.
The computed cross-field $ \mathbf{v}$ is however rather inaccurate.
It is indeed continuous and singularities
are areas in this continuous field where the norm deviation from unity
fails to be ruled out by the penalty term. 
This inaccuracy does not heavily affect the singularity pattern,
because singularity indices are integer quantities and their location
need only to be known approximately. However,
it jeopardizes the accurate subsequent tracing of separatrices
and thus justifies the construction of a second more accurate cross-field
with a truly nonlinear field representation,
and whose computation is based on the singularity pattern
that has been extracted from the first computation
(Sect.~\ref{S:singularities}).

Let the singularity pattern be defined as 
the set $\mathcal{S}=\{S _j, j\in [| 1,N |]\}$ of $N$ singularities located at ${\bf X}_j$ and of respective index $k_j$.
A cross-field  $\mathcal{C}_\Omega$ matching this singularity pattern must fulfill 
three constraints:
\begin{itemize}
\item if $\mathbf{X}$ is on $\pd\Omega$,  at least one branch of $\mathcal{C}_\Omega(\mathbf{X})$ is tangent to $\pd\Omega$,
\item if $\mathbf{X}$ is a singularity of $\mathcal{C}_\Omega$ then $\mathbf{X}\in\mathcal{S}$,
\item if $\mathbf{X}=\mathcal{S}_j$ then $index(\mathcal{C}_\Omega(\mathbf{X}))=k_j$.
\end{itemize}
As in this paper only the {\em planar case} is considered, 
$\Omega \subset \mathbb C$ can be regarded as a subset of the complex plane.
We shall note $\vec n$ the unit normal vector to the complex plane, 
$\vec g$ the interior unit normal vector to the boundary $\pd \Omega$, 
and $\vec T = \vec g \times \vec n$ the unit vector tangent to $\pd \Omega$. 

\subsection{Computation of the  $H$ field}
We are looking for a holomorphic function
\begin{equation}
  \begin{array}{ccrcl}
    \mathcal{F}&:& \mathcal{P}\subset\mathbb{C} & \rightarrow & \mathbb{C}\\
    & &  \vec z=U+\i V & \mapsto & \vec P = x + \i y
  \end{array}
  \label{eq:holomFunc}
\end{equation}
where $\mathcal P$ is a parametric space. As $\mathcal{F}$ and $\mathcal P$ are unknown, solving directly this problem is proved to be a difficult task. Instead, the focus is on finding the $2 \times 2$ jacobian matrix of $\mathcal F$
\begin{equation}
 J_{\mathcal F}(\vec z) = ( \pd_U \mathcal F(\vec z), \pd_V \mathcal F(\vec z) )
 \equiv ( \tilde{\vec u}, \tilde{\vec v} ),
 \label{eq:jacMat}
\end{equation}
where $\tilde{\vec u}, \tilde{\vec v}$ are the column vectors of $J$. 
The 
function $\mathcal F$ being holomorphic, the columns of $J$ have the same norm $|| \tilde{\vec u} || = || \tilde{\vec v} || \equiv e^H$ 
and are orthogonal to each other, $\tilde{\vec u} \cdot \tilde{\vec v} =0$,
so that one can write 
$$
J_{\mathcal F}(\vec z) = e^H  ( \vec u, \vec v )
\quad , \quad 
\vec n = \vec u \land \vec v
$$
$\tilde{\vec u}$ and  $\tilde{\vec v}$ are $2$ branches of the cross-field we are looking for, and $\tilde{\vec u}$ can be represented by a complex function 
\begin{equation}
f:\Omega\rightarrow\mathbb{C}
\quad , \quad 
f = e^{H + \i \theta}.
\label{eq::defBranch}
\end{equation}
The computation of the cross-field is rephrased 
into the computation of two real functions, 
$H : \Omega\rightarrow\mathbb{R}$ and $\theta : \Omega\rightarrow\mathbb{R}$,
from which the first branch $\tilde{\vec u} = f  = e^{H + \i \theta}$ can be deduced, and the second branch of the cross-field is simply $\tilde{\vec v} = \i \tilde{\vec u} = \i f$.

By construction, $\tilde{\vec u}$ being a partial derivative of a holomorphic function, $f$ is holomorphic.
As (\ref{eq::defBranch}) shows, 
$H$ and $\theta$ are the real and imaginary parts of a complex logarithmic function. Indeed,
\begin{equation}
  \begin{array}{lcl}
    \log(f) & = & \ln(|f|) + \i \arg(f)\\
    & = & H + \i (\theta + 2k\pi), k\in\mathbb Z.
  \end{array}
\label{eq::logF}
\end{equation}
The complex logarithm being a holomorphic function, 
both $H$ and $\theta$ are harmonic real functions.
The $H$ field, which is the real part of the logarithm, 
is continuous on $\Omega\setminus \mathcal S$,
i.e., on the domain $\Omega$ from which the pointwise singularities $\mathcal{S}_j$ have been excluded.
The $\theta$ field, on the other hand, is multi-valued
and a {\em branch cut} (see Sect.~\ref{S:cutGraph}) has to be additionally defined on $\Omega\setminus \mathcal S$ 
to have it single-valued, and make its computation possible.
Finally  $H$ and $\theta$ obey the Cauchy-Riemann equations
\begin{equation}
  \nabla \theta = \vec n \times \nabla H,
  \label{eq::H}
\end{equation}
which shall allow to obtain $\theta$, once $H$ has been computed.

The $H$ field is solved first. 
As shown in \cite{Bethuel:1994},
it is the solution of the partial differential problem on $\Omega$
\begin{equation}
  \nabla^2 H = 2 \pi \sum_{j=1}^N \frac{k_j}{4}  \delta ({\bf X}_j).
  \label{eq::HSing}
\end{equation}
The boundary condition $\vec T(s) = (\cos \theta(s), \sin\theta(s))$ ensures that 
the cross, represented by $\theta(s)$,
is tangent to the boundary $\pd \Omega$, represented by the tangent vector $\vec T(s)$,
with $s$ a curvilinear coordinate on $\partial\Omega$. 
This entails by derivation
$$
\pd_s \vec T = (-\sin\theta, \cos\theta) \  \pd_s\theta 
= \vec g \ \pd_s\theta. 
$$
On the other hand, one knows from the differential geometry of surfaces
that
$\pd_s \vec T \equiv \kappa_g \vec g $ 
with $\kappa_g$ the geodesic curvature of the curve $\partial\Omega$. 
One has thus by identification
$$
\kappa_g = \pd_s \theta = \vec T \cdot \nabla\theta 
=  ( \vec n \times \nabla H) \cdot \vec T  = -\vec g \cdot \nabla H,
$$
where Eq.~(\ref{eq::H}) has been used. 

Summing up, the computation of $H$ consists in solving the linear boundary value problem 
\begin{equation}
\left\{\begin{array}{cclll}
\nabla^2 H &=& 2 \pi \sum_{j=1}^N \displaystyle\frac{k_j}{4}  \delta ({\bf X}_j)  & \text{in}& \Omega \\
-\vec g \cdot \nabla H &=& \kappa_g& \text{on}& \partial \Omega.
\end{array}\right.
\label{eq:pbPlanar}
\end{equation}
As the boundary condition is of Neumann type on the entire boundary,
a solution $H$ exists if the condition 
\begin{equation}
  \int_{\pd \Omega} \kappa_g \ \dext\pd\Omega 
  = \int_\Omega  2 \pi \sum_{j=1}^N \frac{k_j}{4}  \delta ({\bf X}_j) \ \dext\Omega 
\label{eq::neumanncompatibility}
\end{equation}
is verified, which is obviously the case,
as the Gauss-Bonnet theorem states the first identity in 
\begin{equation}
  \int_{\pd \Omega} \kappa_g \ \dext\pd\Omega  
  = 2\pi \chi (\Omega) 
  = \int_\Omega  2 \pi \sum_{j=1}^N \frac{k_j}{4}  \delta ({\bf X}_j) \ \dext\Omega
  \label{eq::exPlanar}
\end{equation}
whereas the second identity is the Poincar\'e-Hopf theorem. 

The solution $H$ of the boundary value problem Eq.~(\ref{eq:pbPlanar}) 
is however known up to an arbitrary additive constant,
which is not harmful as only $\nabla H$ is needed to determine $\theta$. 
It is a linear problem solved using the finite element method 
on a triangulation of
$\Omega$ with order 1 Lagrange elements.
Once $H$ is known, one proceeds with the determination of $\theta$,
as explained in the next section, to complete the cross-field $(f,\imath f)$. 

\begin{figure}[h!t]
\begin{center}
  \includegraphics[width=0.35\textwidth]{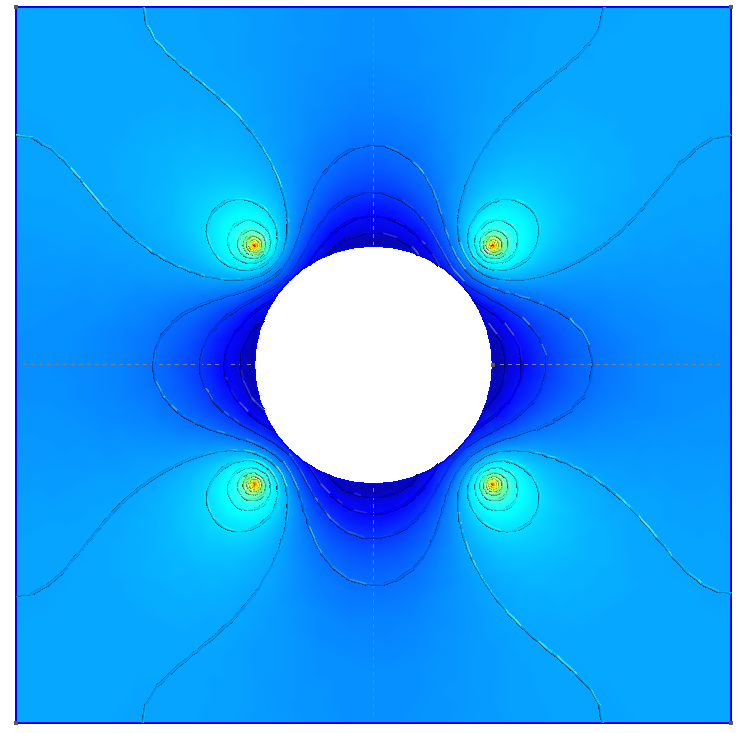}
\caption{$H$ field obtained for 4 singularities of valence $3$ on corners, 
  and 4 singularities of valence $5$ inside the domain.}
\label{fig:exH}
\end{center}
\end{figure}

\subsection{Computation of the $\theta$ field}
\label{S:Hcomp}

As the curl of the right-hand side of Eq.~(\ref{eq::H}) 
$$
\nabla \times (\vec n \times \nabla H) = (\nabla^2 H) \vec n = \vec 0
$$ 
vanishes in $\Omega\setminus \mathcal S$, 
one knows that the left-hand side is indeed a gradient,
and that there exists a scalar field $\theta$ 
verifying the differential equation (\ref{eq::H}).
This function is however multi-valued
and a branch cut has to be determined to make it single-valued, and make its computation possible.

\subsubsection{Generate the branch cut}
\label{S:cutGraph}
A branch cut is a set $\mathcal{L}$ of curves of a domain $\Omega$ that do not form any closed loop 
and that cut the domain in such a way that it is impossible to find any closed loop in $\Omega\setminus \mathcal{L}$ that encloses one or several singularities, or an internal boundary. 
As we already have a triangulation of $\Omega$, the branch cut $\mathcal{L}$ 
is in practice simply a set of edges of the triangulation. 
The method presented here is based on \cite{Bommes:2009}.

\begin{figure}[h!t]
\begin{center}
  \includegraphics[width=0.35\textwidth]{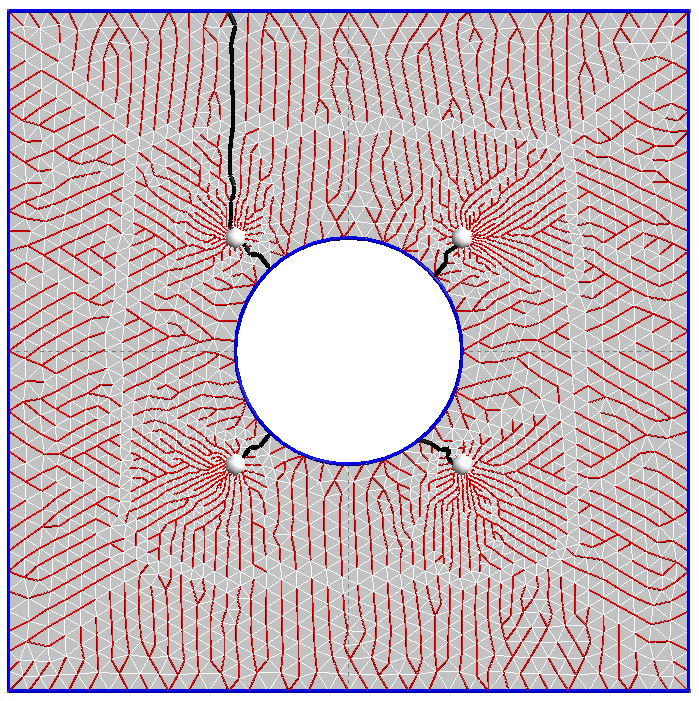}
\caption{$\pd \Omega$ is represented in blue, in white the singularities, in red are edges of the spanning tree $\mathcal{T}$, and in black the edges of the branch cut $\mathcal{L}$.}
\label{fig:cutGraph}
\end{center}
\end{figure}

The nodes and the edges of a triangulation of a domain $\Omega$ can be regarded as a graph $\mathcal{G}$.
A spanning tree $\mathcal{T}$ is a sub-graph of $\mathcal{G}$ that links all vertices of $\mathcal{G}$ with the minimum number of edges.
A spanning tree can be generated by a Breadth-first-search algorithm. One starts from an arbitrary vertex $v$, and adds to the spanning tree $\mathcal{T}$ all edges $(v, v_i)$ that are adjacent to $v$ and that connect vertices $v_i$ that have not yet been visited by the algorithm. The procedure is continued recursively with all newly considered vertices $v_i$ until running out of unvisited vertices. 
There exist many equivalent spanning trees for a given mesh. Their shape can be improved by modifying the Breadth-first-search algorithm so that vertices which are the closest to the singularities or to the boundaries are processed in priority. A spanning tree built this way is presented in red in Fig.~\ref{fig:cutGraph}.

This spanning tree forbids however closed curves around all nodes of the mesh, which is an unnecessarily strong restriction,
as we only want to forbid closed curves around the singularities. 
If we call hanging edge of $\mathcal{T}$ an edge whose leaf node $v_l$ is not a singularity
(a leaf node is a node of the graph adjacent to exactly one edge), 
the branch cut is the sub-graph of $\mathcal{T}$ obtained by successively removing hanging edges until no hanging edges is left in the spanning tree. The result of this substraction is the branch cut $\mathcal{L}$ needed for the computation of $\theta$.
It is depicted in black in Fig.~\ref{fig:cutGraph}. 

\subsubsection{Solve the Cauchy-Riemann equations}
\label{S:theta}

Once a branch cut $\mathcal{L}$ is available, the field $\theta$ can be computed by solving the linear Cauchy-Riemann equations (\ref{eq::H}). The chosen boundary condition consists in fixing the angle $\theta$ at {\em one} arbitrary point $\vec P \in \pd \Omega$ so that  $\mathcal{C}_\Omega (\vec P)$ has one of its branch collinear with $\vec T (\vec P)$.
The problem can be rewritten as the well-posed Eq.~(\ref{eq:thetaPlanar}) and is solved using finite element method on the triangulation $\Omega_T$ with order one Crouzeix-Raviart elements. This kind of elements have proven to be more efficient for cross-field representation as underlined in \cite{Jezdimirovic:2019}.

\begin{equation}
  \left\{\begin{array}{rl}
  \nabla \theta & = \vec n \times \nabla H \text{ in } \Omega \\
  \theta(\vec P) & = \theta_{\vec P} \text{ for an arbitrary } \vec P \in\partial \Omega\\
  \theta & \text{ discontinuous on } \mathcal{L}
\end{array}\right.
\label{eq:thetaPlanar}
\end{equation} 

\begin{figure*}[h!t]
\begin{center}
  \includegraphics[width=0.32\textwidth]{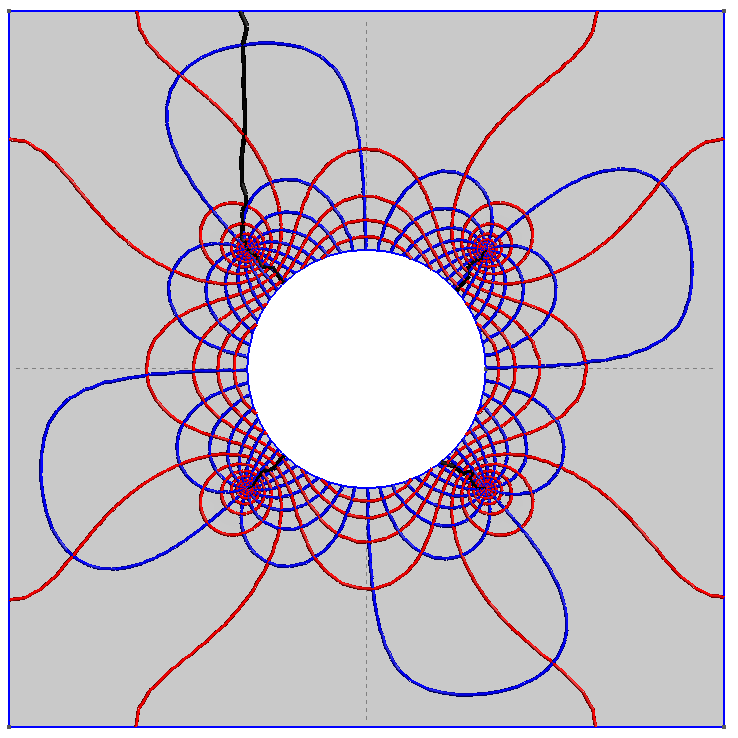}
  \includegraphics[width=0.32\textwidth]{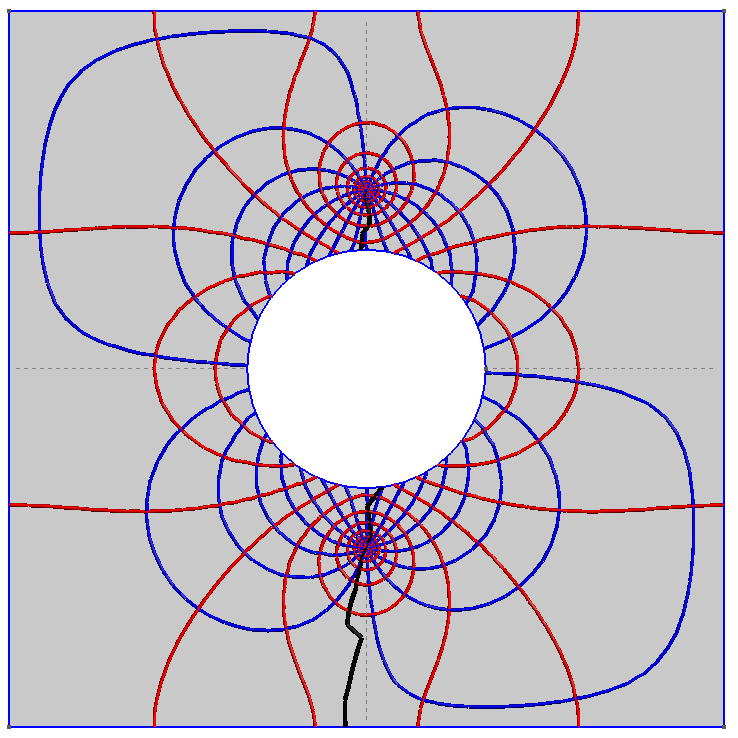}
  \includegraphics[width=0.32\textwidth]{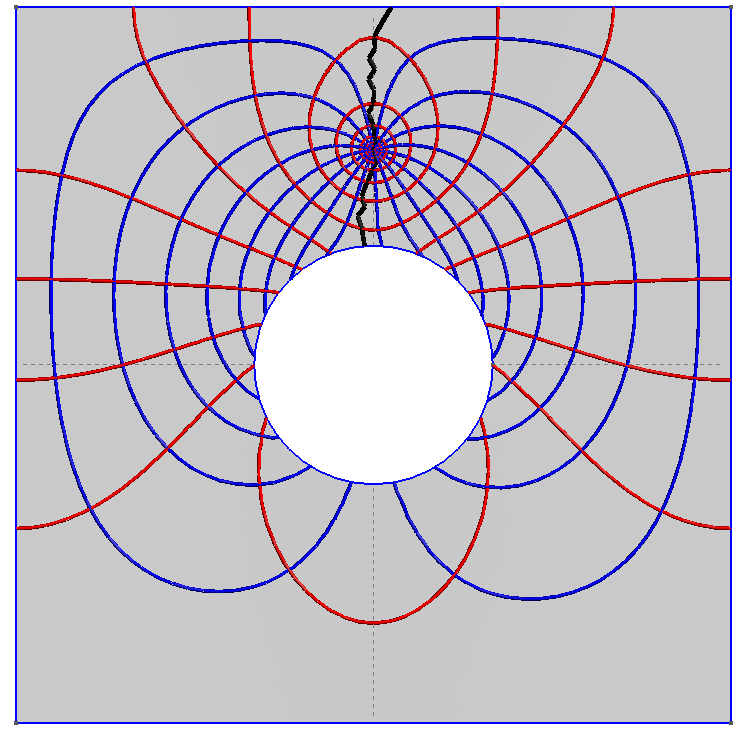}
\caption{$H$ isolvalues (in red) and $\theta$ isovalues (in blue) for different singularities configurations. The black line represent the cutgraph. All cases have valence $3$ singularities imposed at corners (for a total of $4$ singularities in each case). The extra singularities which are added inside the domain are: $4$ singularities of valence $5$ (left figure),  $2$ singularities of valence $6$ (center figure),  $1$ singularity of valence $8$ (right figure).}
\label{fig:thetaH}
\end{center}
\end{figure*}

Fig.~\ref{fig:thetaH} presents results for $3$ different types of singularities. As expected, $H$ and $\theta$ isolines are orthogonal, this two scalar fields being harmonic conjugate.

Once $H$ and $\theta$ scalar fields are computed on $\Omega$, the cross-field $\mathcal{C}_\Omega$ can be retrieved for all $\vec X \in \Omega$:

\begin{equation}
\vec c (\vec X) = \{ \vec u_k = e^{H(\vec X)+jk\theta(\vec X)}, k\in[|1,4|]\}
\label{eq:HThetaToCross}
\end{equation}

\subsubsection{Non quad-meshable cross-fields}
\label{S:nonMeshable}

We showed in the previous section that if the imposed singularity pattern respects the condition on Euler's characteristic, then there exists a unique scalar function $H$ (up to an additive constant) verifying Eq.~(\ref{eq:pbPlanar}) and a unique scalar function $\theta$ verifying Eq.~(\ref{eq:thetaPlanar}). This is leading us to a unique cross-field solution $\mathcal{C}_\Omega$ to our problem.

However, we never showed that the obtained cross-field $\mathcal{C}_\Omega$ is actually respecting the following and mandatory condition:

\begin{itemize}
\item $\forall \mathbf{X}\in\pd\Omega$,  at least one branch of $\mathcal{C}_\Omega(\mathbf{X})$ is tangent to $\pd\Omega$
\end{itemize}

Indeed, problem (\ref{eq:thetaPlanar}) is well-posed and only ensures this condition for a given $\vec P \in \pd\Omega$. Adding extra boundary conditions for $\theta$ will lead to make problem (\ref{eq:thetaPlanar}) ill-posed.Therefore, when domain $\Omega$ is not simply connected, it is possible that obtained a cross-field $\mathcal{C}_\Omega$ not tangent on parts of $\pd\Omega$.

\begin{figure}[H]
\begin{center}
  \includegraphics[width=0.3\textwidth]{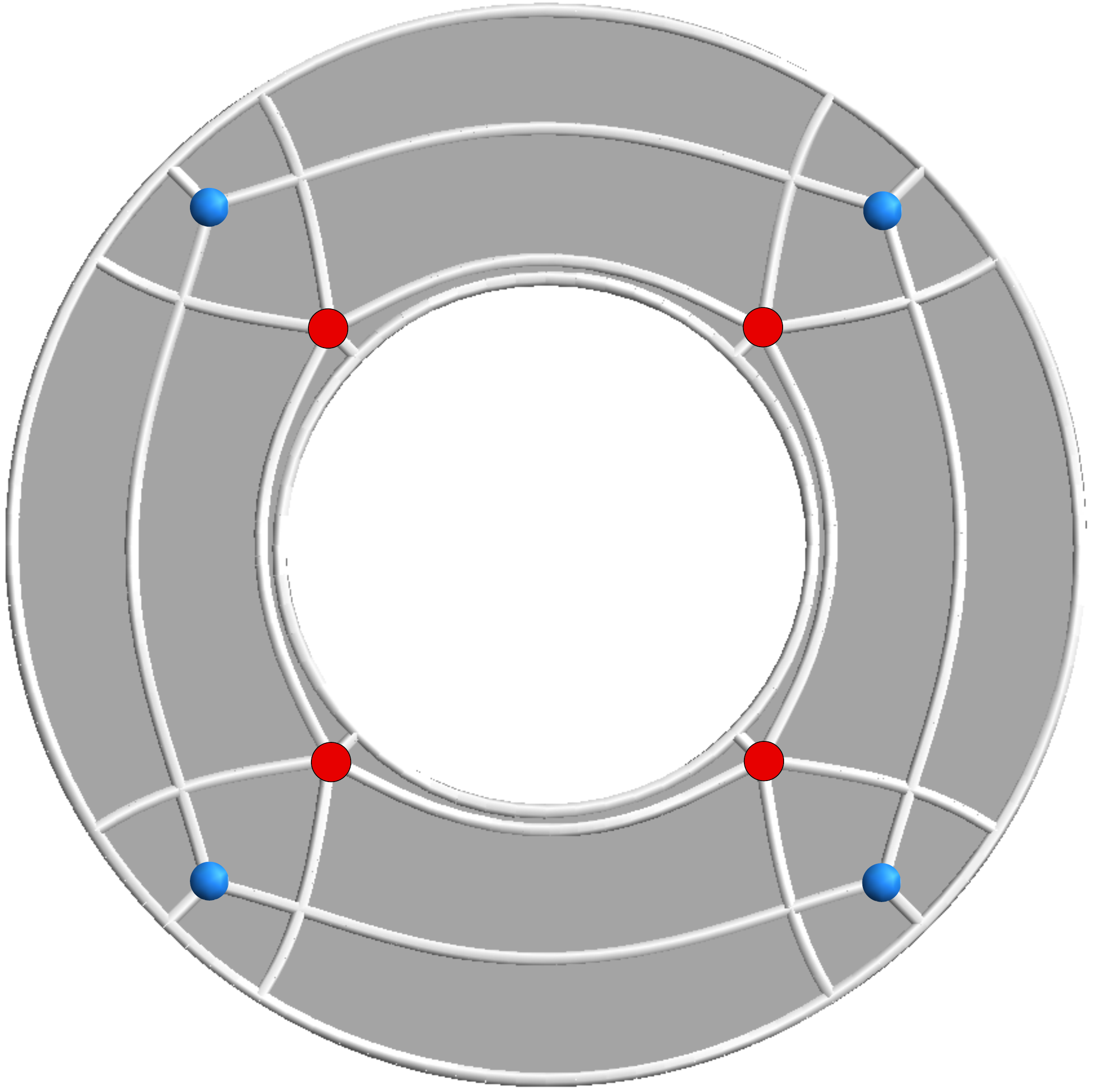}
  \includegraphics[width=0.3\textwidth]{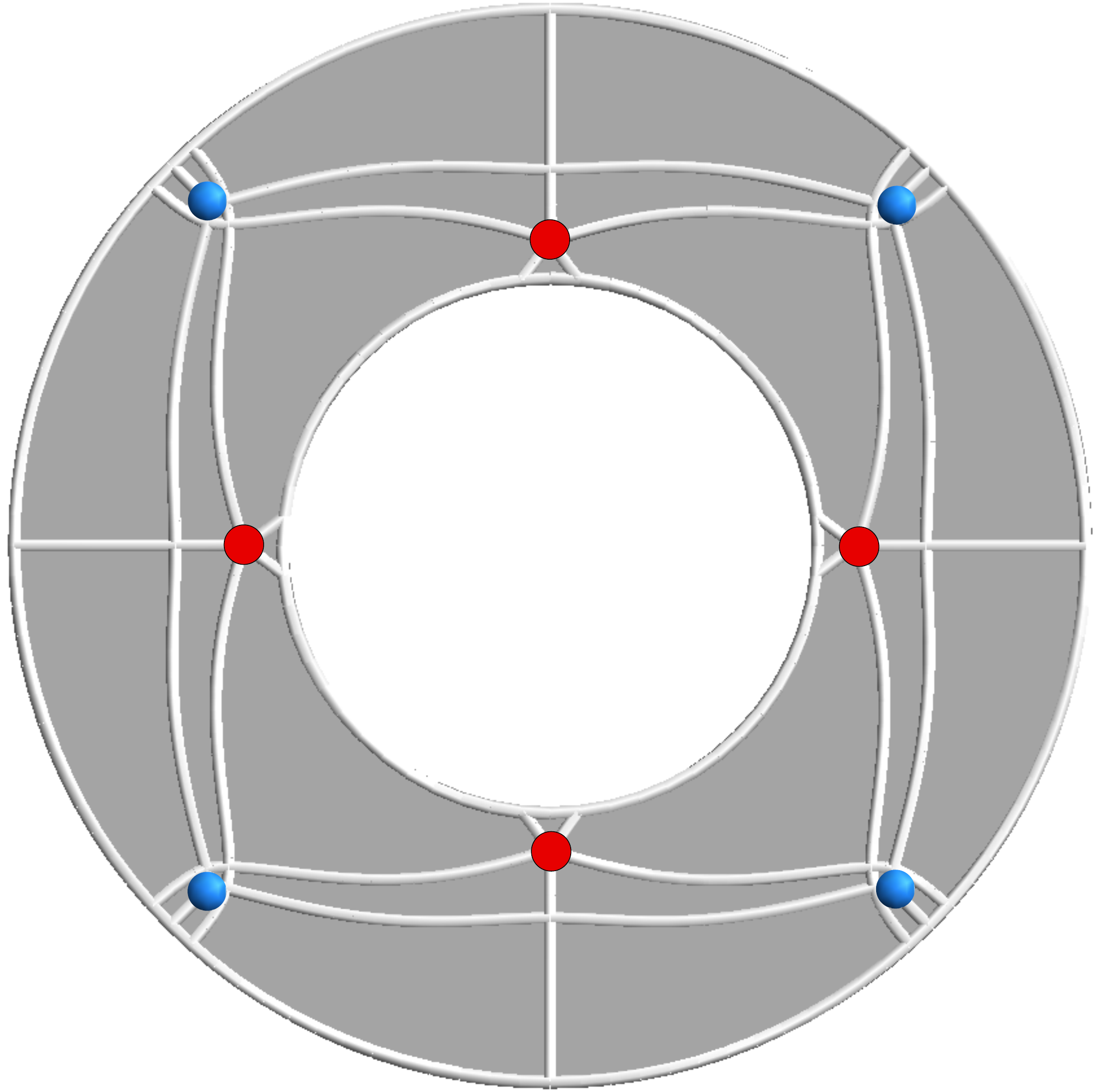}
\caption{Highlighting of singularity patterns for which no holomorphic transformation matching all imposed conditions exists.}
\label{fig:nonMeshable2}
\end{center}
\end{figure}

When this happens, it means that there is no holomorphic function $\mathcal F$ (eq.~\ref{eq:holomFunc}) matching all conditions imposed ({\em i.e.} singularity pattern and tangency to domain boundaries). Figure~\ref{fig:nonMeshable2} shows two different singularity patterns where only the singularity location differs for a same domain $\Omega$. The domain $\Omega$ considered is an annulus, and a set of $4$ index $-1$ and $4$ index $1$ singularities is imposed. Singularity pattern for the example on the left leads to a quad decomposition of the domain and singularity pattern for the example on the right does not. Cross-field obtained for example on the right does not respect the tangency condition on $\Omega$ interior boundary. There exists no holomorphic transformation matching at the same time singularity pattern and cross-field tangency to $\Omega$ boundaries.

For the issue of imposing correct singularity configuration we refer the user to Abel-Jacobi conditions given in \cite{Lei:2020}.

\subsection{Computing the per-partition parameterization}
\label{subsec:parametricSpace}

Obtaining a global parameterization comes down to computing real values functions $U(x,y)$ and $V(x,y)$ defined in Eq.~(\ref{eq:holomFunc}), from the known smooth jacobian of the mapping $\mathcal F:$ $J_{\mathcal F}(\vec z) = ( \pd_U \mathcal F(\vec z), \pd_V \mathcal F(\vec z) ) \equiv ( \tilde{\vec u}, \tilde{\vec v} )$ Eq.~(\ref{eq:jacMat}). Establishing $U$ and $V$ is equivalent to determining $\mathcal F ^{-1}$. We know that $\forall \vec P = \mathcal{F}(\vec z) \in\Omega\setminus \mathcal S$, $J_{\mathcal F}$ corresponds to two branches of the cross-field.

We also know that $J_{\mathcal F^{-1}}$ jacobian of $\mathcal{F}^{-1}$ is, $\forall \vec P  \in \Omega\setminus \mathcal S$, such as:
\begin{equation}
  \left\{
  \begin{array}{lcl}
  J_{\mathcal F^{-1}}(\vec P)& = & (\pd_x \mathcal F^{-1}(\vec P), \pd_y \mathcal F^{-1}(\vec P)) \\
    & = & J^{-1}(\mathcal{F}^{-1}(\vec P)) \\
   &  \equiv & (\bar{\vec u}, \bar{\vec v} )^T.
   \end{array}    \right.
  \label{eq:invJ}
\end{equation}

which, by using the Eq.~(\ref{eq:jacMat}), gives:
\begin{equation}
  \left\{
  \begin{array}{lcl}
    \bar{\vec u} &= e^{-H} \vec{u} &= \nabla U\\
    \bar{\vec v} &= e^{-H} \vec{v} &= \nabla V.\\
  \end{array}
  \label{eq:gradInv}
  \right.
\end{equation}

Now, the parameterization can be simply obtained by solving:
\begin{equation}
\begin{cases}
\nabla U = \bar{\vec u} \\
\nabla V = \bar{\vec v},
\end{cases}
\label{eq:paramUV}
\end{equation}
using a finite element method with order one Lagrange elements.

We mentioned earlier that for the problem to be well posed, $J_{\mathcal{F}}$ has to be smooth, meaning that $\tilde{\vec u}$ and $\tilde{\vec v}$ have to be smooth. This implies that a lifting of the cross-field $\mathcal C_\Omega$ has to be done on $\Omega$ allowing discontinuities of $(\tilde{\vec u},\tilde{\vec v})$ across $\mathcal L$. Once this operation is done, it is possible to determine $(\bar{\vec u},\bar{\vec v})$ and compute $U$ and $V$ by solving the Eq.~(\ref{eq:paramUV}). By construction, $U$ and $V$ will be discontinuous across $\mathcal L$, and their isolines will be in all points tangent to the cross-field $\mathcal C_\Omega$.

\begin{figure}[h!t]
\begin{center}
  \includegraphics[width=0.32\textwidth]{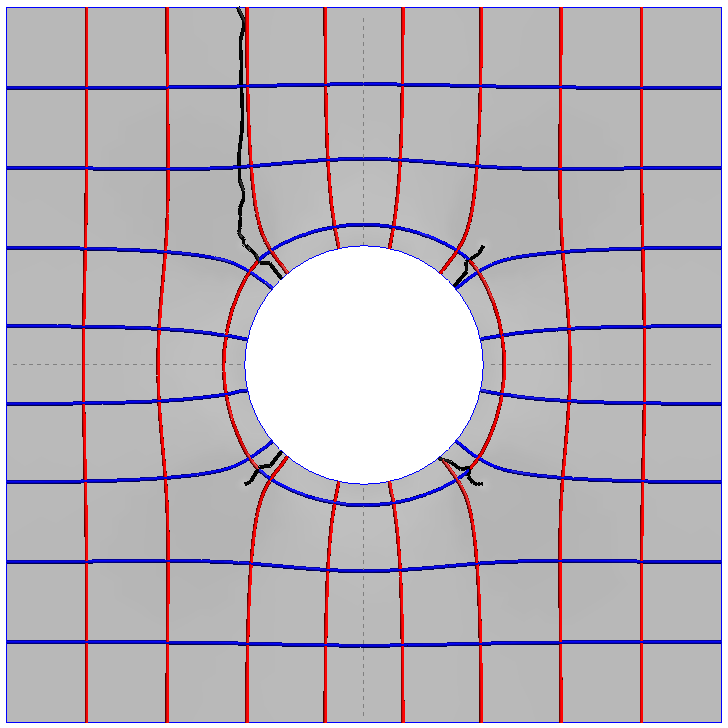}
  \includegraphics[width=0.32\textwidth]{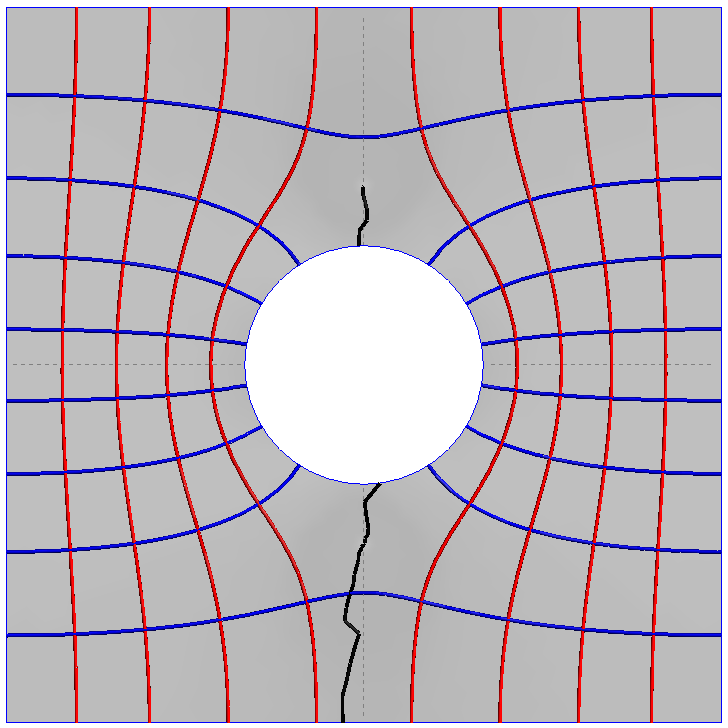}
  \includegraphics[width=0.32\textwidth]{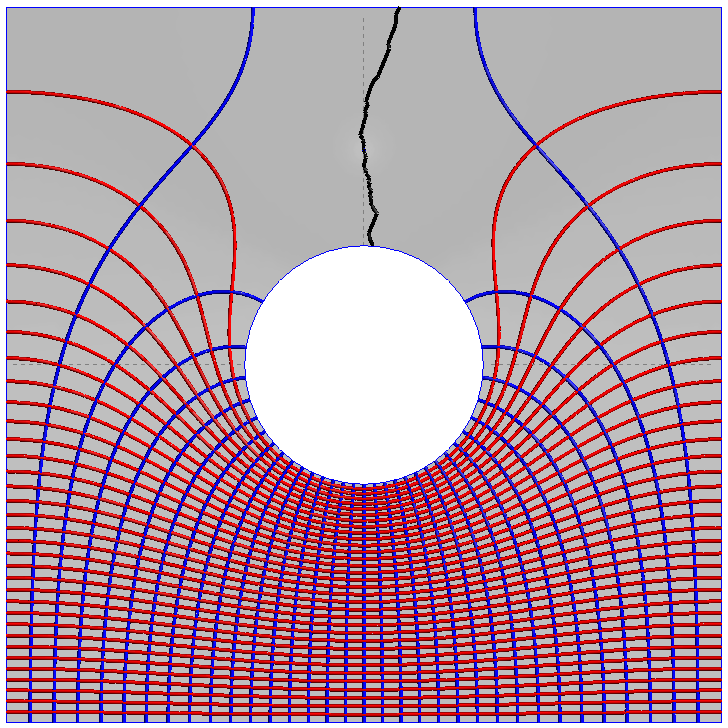}
\caption{$U$ isolvalues (in red) and $V$ isovalues (in blue) for different singularities configurations. The black line represents the cutgraph.}
\label{fig:UViso}
\end{center}
\end{figure}

It is important to note that for this global parameterization computation, no conditions are imposed on the cut graph $\mathcal L$. Therefore, nothing guarantees correspondence of $U$ and $V$ isolines across $\mathcal L$. 
In principle, this global parameterization is not enough by itself to generate a conforming quadrilateral mesh of a domain $\Omega$. Some works 
bypass these difficulties
\cite{Bommes:2013, Bommes:2009, Kalberer:2007} by constraining parametric coordinates on boundaries and on $\mathcal L$, however it is not the direction followed in our approach.

\begin{center}
\begin{figure}[!htb]
\begin{center}
\includegraphics[width=0.7\textwidth]{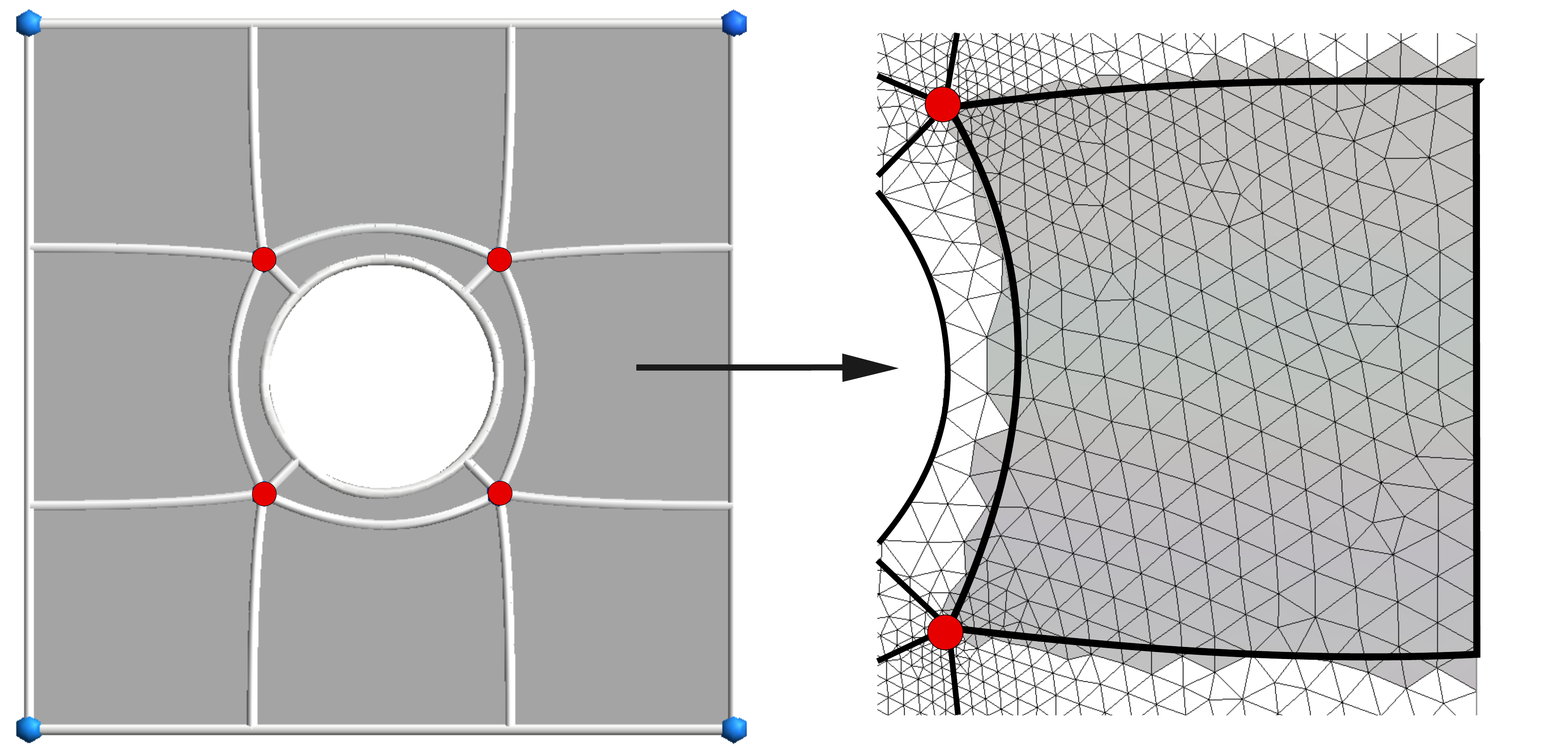} 
\caption{Extracted triangulated partition for parameterization}
\label{fig:triEdges}
\end{center}
\end{figure}
\end{center}

The goal here is to obtain per-partition parameterization, starting from the construction of a conformal quad layout
 of $\Omega$ (section \ref{S:MultiBlock}) and followed by the computation of independent parameterizations for each of the resulting quadrilateral partitions. The per-partition parameterization computed in this manner is aligned with smooth cross-field (singularities can only be located on corners of the partitions). To achieve this goal, the first step is to extract corresponding triangles of each patch. Due to the mutual separatrices' points, neighboring partitions contain the overlapping triangles, as shown in Fig.~\ref{fig:triEdges}.

\begin{center}
\begin{figure}[!htb]
\begin{center}
\includegraphics[width=0.7\textwidth]{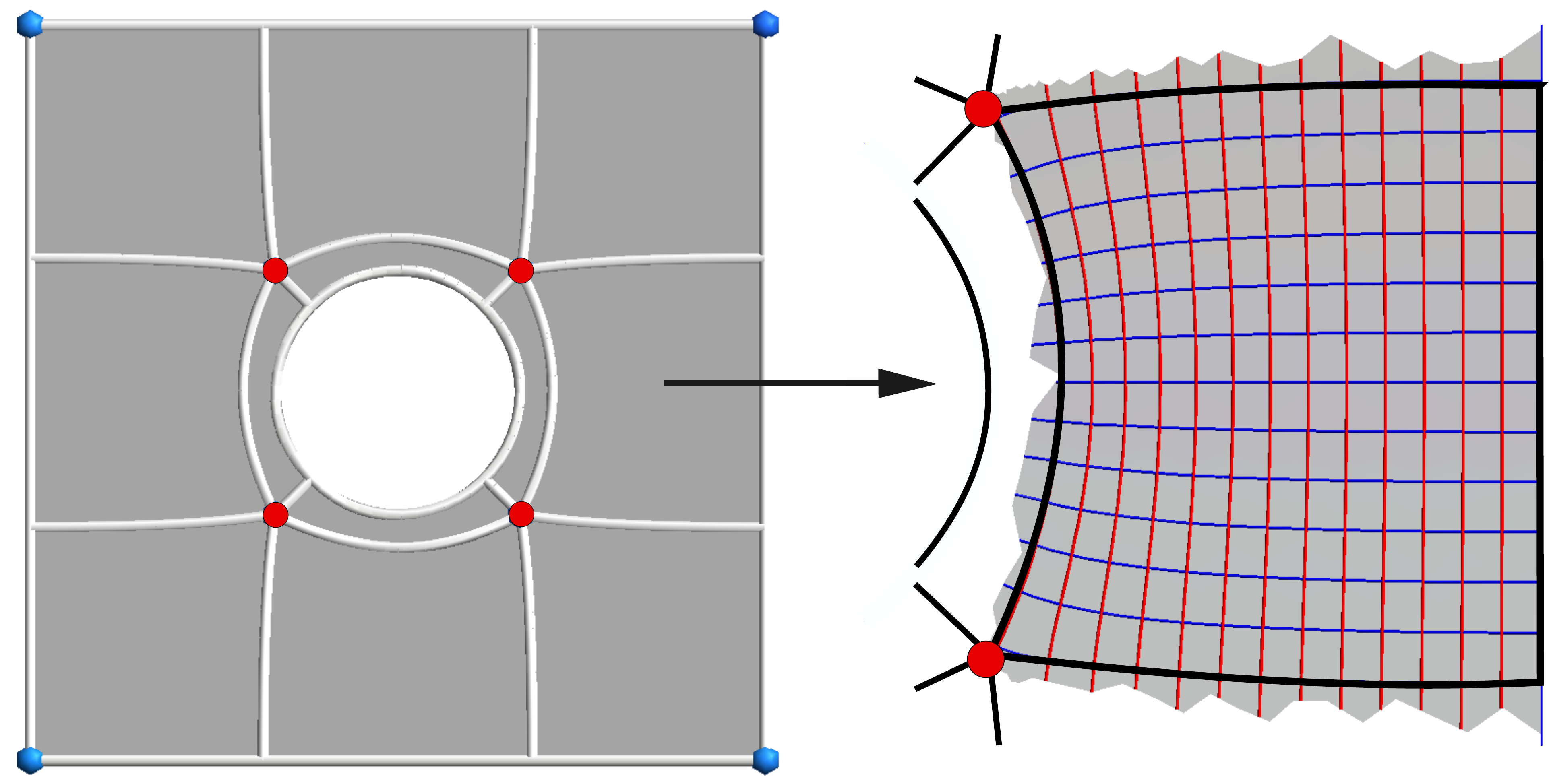} 
\caption{$U$ isolvalues (in red) and $V$ isovalues (in blue) on extracted partition}
\label{fig:parametrization}
\end{center}
\end{figure}
\end{center}

An important remark is that no cut graph is needed here to perform a lifting of the underlying cross-field to obtain smooth couple of vector fields $(\tilde{\vec u},\tilde{\vec v})$. As no cut graph is used, all isolines of $U$ and $V$ will be continuous inside the patch, Fig~\ref{fig:parametrization}.

Based on \cite{Jezdimirovic:2019} for locally integrable cross-fields, note that on each partition holds:
\begin{equation}
\begin{cases}
\nabla \times e^{-H}u=0 \\
\nabla \times e^{-H}v=0,
\end{cases}
\end{equation}
where from $(\tilde{u}, \tilde{v})$ is a locally curl-free cross-field. Thus, the resulting parameterization $(U, V)$ is \emph{conformal} on each patch (\cite{Ray:2006}). Following \cite{Floater:2005}, a  \emph{conformal} mapping is also \emph{harmonic}. Radó-Kneser-Choquet (RKC) theorem for planar harmonic mappings (details in \cite{Floater:2005, Iwaniec:2019}) states:

If $f : \mathcal{S} \rightarrow \mathbb{R}^2$ is harmonic and maps the boundary $\partial \mathcal{S}$ homeomorphically into the boundary $\partial \mathcal{Q}$ of some convex region $\mathcal{Q} \subset \mathbb{R}^2$, then $f$ is \emph{one-to-one}.

In our case, the parametric (image) region is in the shape of a rectangle, thus, our parameterization $(U, V)$ is per-partition \emph{bijective}.

\section{Quad layout: generating and correcting the partitions}
\label{S:MultiBlock}

In order to reach the ultimate goal of generating a block-structured quad mesh, we must ensure that the quad layout contains only conformal quadrilateral partitions. To do so, a special attention is paid on obtaining the minimal number of separatrices, existence of limit cycles and non-quadrilateral partitions. In the following, we present ways to modify all ill-defined partitions, such that the overall quad layout topology is not jeopardized.

\subsection{Obtaining the partitions}
\label{MBgenerating}
To obtain the quad layout, the separatrix-tracing algorithm is performed on retrieved cross-field from $H$ function. Initiation of separatrices in singular areas follows the method from \cite{Jezdimirovic:2019}, while the tracing of separatrices in the smooth parts of domain relies on Heun's (a variation of Runge-Kutta 2) numerical scheme.

\begin{figure}[h!t]
\begin{center}
\includegraphics[width=0.55\textwidth]{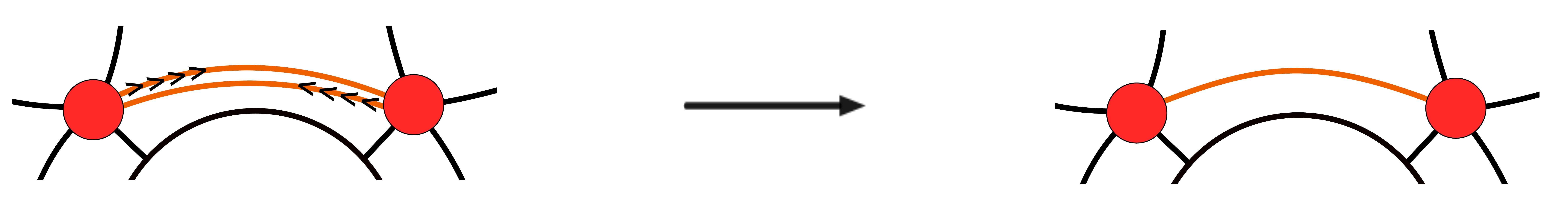} 
\caption{Illustration of doubled separatrice and removal of the duplicate}
\label{fig:doubleSep}
\end{center}
\end{figure}

Our algorithm doesn't have an order of tracing the separatrices - they are all generated simultaneously and independently, similarly to \textit{motorcycle graph} algorithm. Contrary to \textit{motorcycle graph's} stopping criteria, in our case, the tracing of separatrices is finished as soon as they reach the boundary or encounter a singular point. The number of generated separatrices in this manner is not appropriate for a quad layout, due to the  separatrices' initiation process which allows tracing the same separatrice from two different singular points, as shown in Fig.~\ref{fig:doubleSep}. To obtain the valid quad layout, criteria for gaining the minimal number of separatrices are adopted from \cite{Jezdimirovic:2019}. 

\begin{figure}[h!t]
\begin{center}
\includegraphics[width=0.55\textwidth]{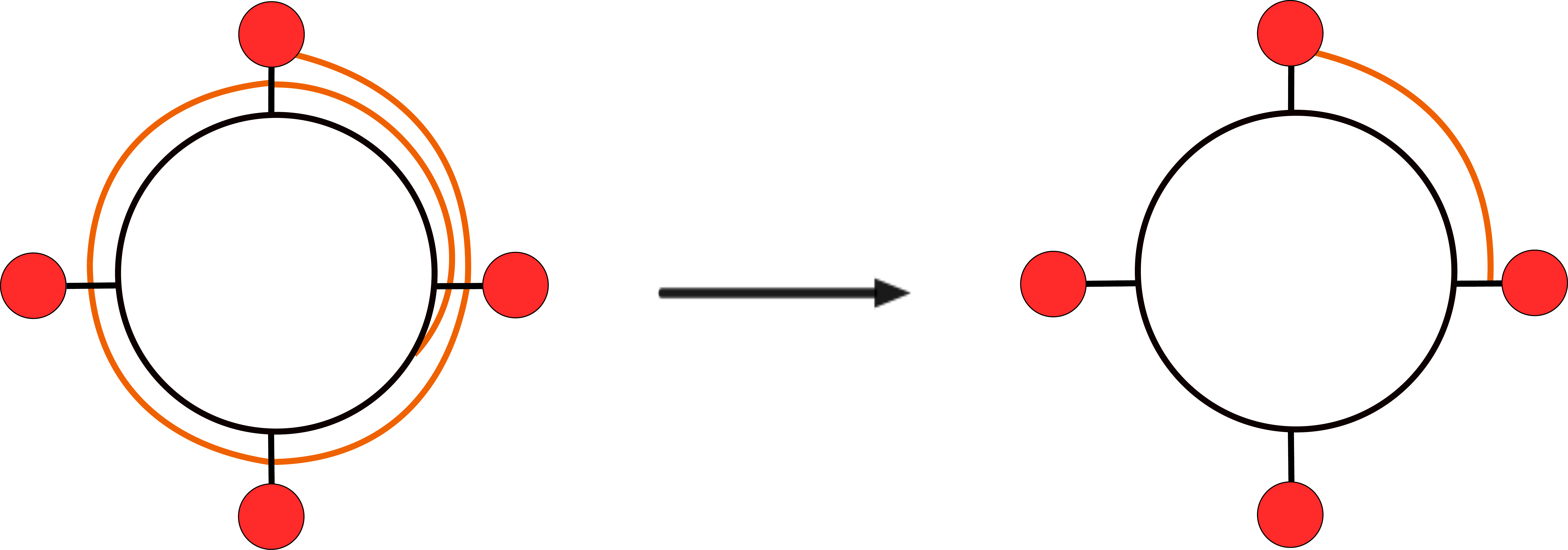} 
\caption{Illustration of an authentic limit cycle treatment}
\label{fig:LimitCycle}
\end{center}
\end{figure}

The next step is to establish the existence of limit cycle(s) in the set of the generated separatrices. In more details, we define a separatrix which is transpassing the other separatrices more than once as a \textit{possible} limit cycle. An \textit{authentic} limit cycle is a possible limit cycle which has more far-flung intersections than its pair, or a separatrix not reaching the boundary orthogonaly. An authentic limit cycle is then cut at the closest orthogonal intersection with another separatrix, as shown in Fig.~\ref{fig:LimitCycle}. By doing so, simplifying the quad layout is accomplished (due to the disappearance of the thin chords), however the obtained partitions can contain T-junctions, which will be modified later on. 

\subsection{Valid quad layout}

If the quad layout contains T-junctions (generated from cutting the authentic limit cycles) or non-quadrilateral partitions, its no longer suitable for block-structured quad meshing.
More specifically, a high-quality quad mesh, according to \cite{Bommes:2012}, must obtain not only the adequate configuration of singularities (in terms of their number and placement) but also respect the connectivity constraint
of being conforming and purely quadrilateral, which the previous issues hinder.

The appearance of T-junctions in our algorithm has already been elucidated, which is not the case with the non-quadrilateral blocks though. In order to examine (and later on correct) the latter case, it's important to mention that our algorithm deals with planar geometries on which, as stressed out by \cite{Bommes:2012} and \cite{Ray:2008}, degenerate singularities (of valence $1$ - \textit{singlet} or of valence $2$ - \textit{doublet}) may appear. For instance, the typical placement of singularities of valence $2$ is on corners of convex boundaries (\cite{Bommes:2012}) or on very sharp/obtuse angles due to the number of cross-field's rotations around the local normal, as the exact formulas for singularity index in \cite{Viertel:2019} show. While these types of singularities do not endanger topological validity of the obtain quad layout, they can cause problems for the block-structured quad meshing by forming non-quadrilateral blocks, as shown in Fig.~\ref{fig:triCorr1}.

\subsection{Correcting the partitions}
\label{S:correction}

\begin{figure}[h!t]
\begin{center}
\includegraphics[width=0.8\textwidth]{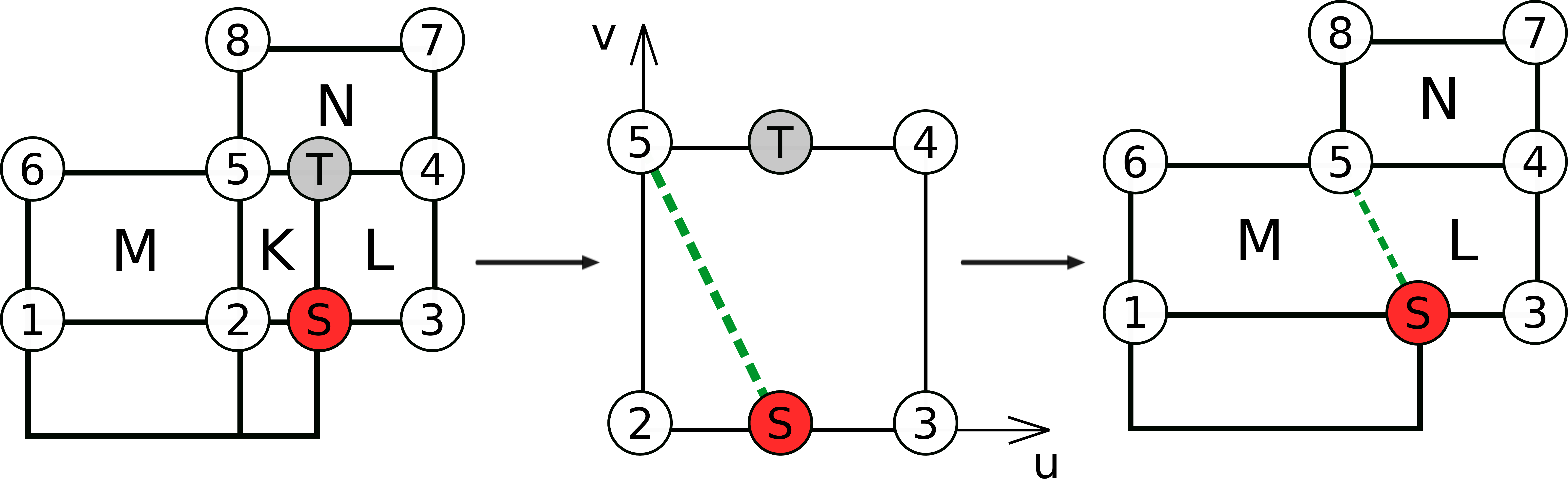} 
\caption{Modifying patches containing T-junctions and affected partitions}
\label{fig:tJunctionCorr}
\end{center}
\end{figure}

When the partition contains a T-junction, the idea is to find the closest neighbor  (which is not the boundary vertex or another T-junction) and, by knowing the parameterization of these blocks, to merge these two separatrices in the parametric space and correct all patches affected by this action, Fig~\ref{fig:tJunctionCorr}. This process is iteratively repeated until the quad layout is T-junction free. 
Following the definition of separatrices and singular points, each of the newly obtained partitions has smooth cross-field inside allowing further on parameterization.

\begin{center}
\begin{figure}[h!t]
\begin{center}
\includegraphics[width=0.75\textwidth]{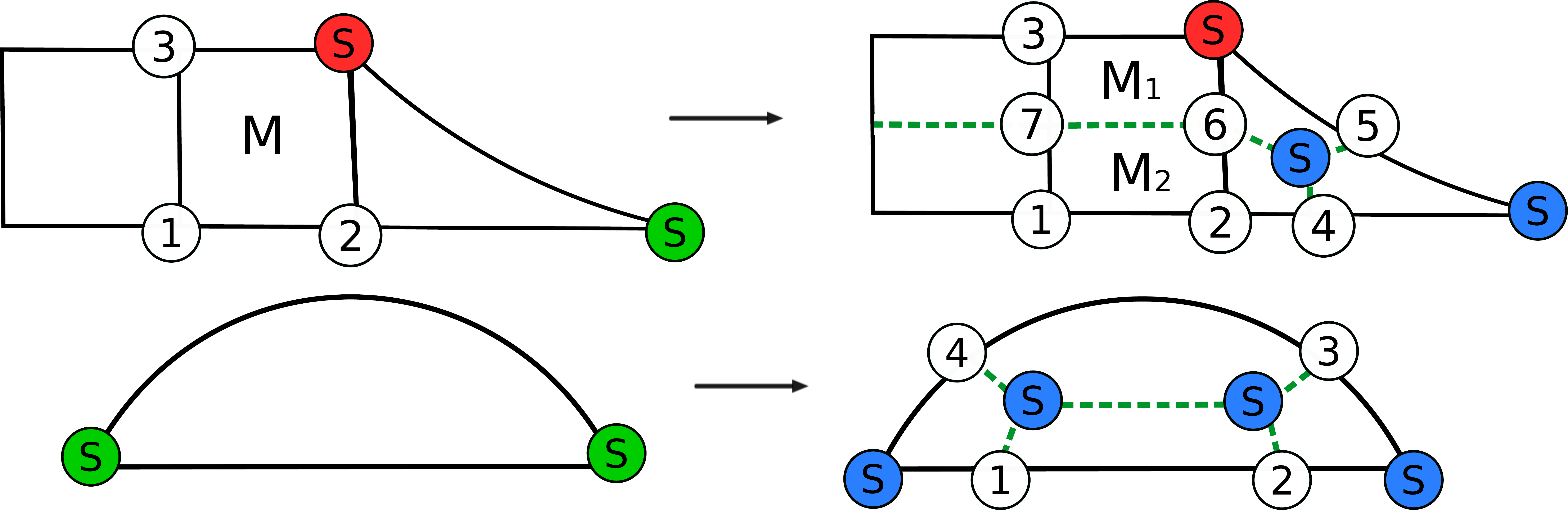} 
\caption{Modifying the patch containing a singularity of valence 2 and affected partitions}
\label{fig:triCorr1}
\end{center}
\end{figure}
\end{center}

Correcting the triangular partitions is performed by forcing the singularity of valence $2$ on the boundary into splitting in two singularities of valence $3$, so that the topology of the overall obtained partitions remains valid (Fig.~\ref{fig:triCorr1}). The position of the imposed singularity is determined simply by finding the barycenter of the partition. Final quad layout is obtained by splitting the triangular patch into three quadrilateral partitions and correcting all patches influenced by the former action.

\newpage
\section{Quad meshing}
\label{S:QuadMesh}

To generate a final quad mesh of the model, we start by remeshing each layout's partition with bilinear transfinite interpolation (TFI) in the parametric space, while respecting the discretizations of layout's edges implied by the field $H$. This step is followed by a computation of a per-partition parameterization aligned with the cross-field (as detailed in section~\ref{subsec:parametricSpace}). Finally, generating a quad mesh of pleasing  quality and intact singularities' positions (examples in section \ref{S:Conclude}) is obtained by mapping TFI mesh onto the physical space.
At this point, it is indeed possible to further optimize vertices' placements via smoothing (e.g. Winslow from \cite{Winslow:1966,Knupp:1999}), in the cases when the application allows/encourages moving singularities away from their initial positions. For the sake of completeness, details of the above-mentioned steps are elucidated in the following.

In order to remesh each patch of the quad layout, the chords' data has been extracted. Keeping in mind that the quad layout is now composed of only conforming quadrilaterals, it is possible to  consistently discretize all its edges by following the size map implied by H:

\begin{equation}
s = e^{H},
\label{eq:sizeMap}
\end{equation}  

where $s$ represents the mesh size. This step tends to generate a uniform quad mesh, unless a specific size map (via implying an adequate configuration of singularities) is imposed by the user. With obtained discretization of edges, a quad mesh in parametric space is trivially created by remeshing each partition with bilinear transfinite interpolation (\cite{Thompson1998}). 

As the parameterization is now defined, the mapping of the quad mesh from parametric onto physical space is achieved by simply finding the triangle $T_M$ in the parametric space to whom the point 
$M(U, V)$ belongs. Afterwards, the computing of coordinates $M(x, y, z)$ is done by assuming that there exists a \emph{linear mapping} from physical space $X$ to parametric space $P$ inside the triangle $T_M$.

\section{Conclusion and Future Work}
\label{S:Conclude}

\begin{figure*}[!htb]
\begin{center}
\includegraphics[width=1.0\textwidth]{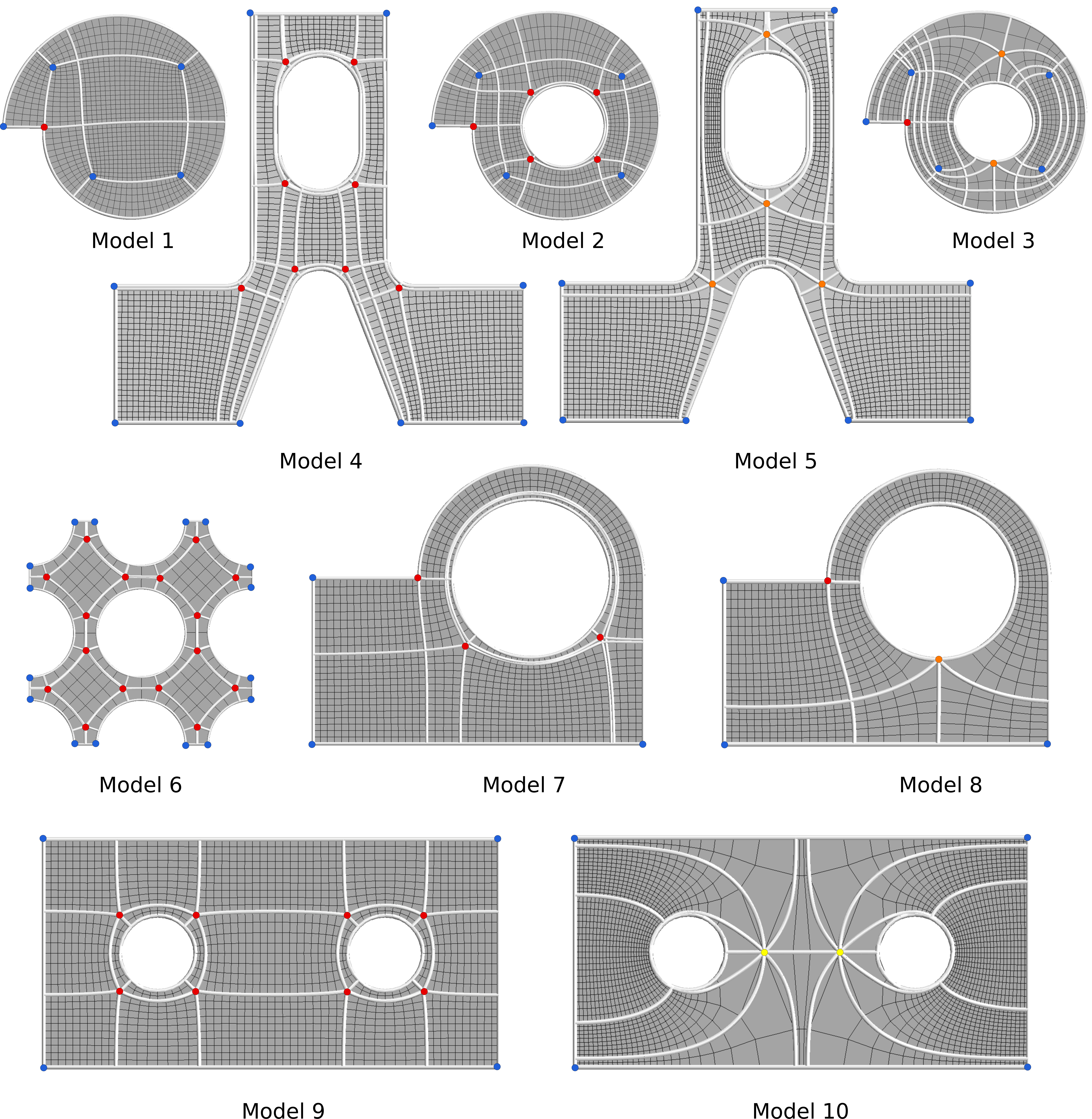} 
\caption{Examples of quad meshes generated from natural and user-imposed singularity configurations. Singularities of valence $3$ are marked in blue, of valance $5$ in red, of valance $6$ in orange and of valence $8$ in yellow.}
\label{fig:meshes}
\end{center}
\end{figure*}

We presented an algorithm that can effectively generate a conformal quad layout from naturally appearing as well as user-imposed set of singularities, as shown in Fig.~\ref{fig:meshes}. The latter case is possible under constraints of the topological correctness and existence of a conformal mapping for imposed singularity pattern. In particular, the developed method is relying on the computation of the two cross-fields in order to accommodate high gradients and improve cross-field representation around singularities. With the correction scheme for limit cycles and non-quadrilateral partitions, we are able to produce conformal quad layout and finally a block-structured quad mesh. Remeshing of each patch is relying on construction of a per-partition parameterization  and following the isolines from the second cross-field. A quad mesh obtained in this manner preserves singularity placement and is of pleasing quality, as shown in Table~\ref{table_meshQuality}. At last, our novel pipeline exhibits interactive design while being simple, automatic and available in the Gmsh software.

Attractive directions for future work include, although are not limited to, extending the algorithm for computations on manifold surfaces, which is currently underway, as well as its adaptation for working on varying user-prescribed element sizes.

\begin{table}[h!]
\begin{center}
\begin{tabular}{lccccc}
\hline
\textbf{Geometry} & \textbf{Edge length} &  \multicolumn{3}{c}{\textbf{Mesh quality}} \\ && $\overline{\eta}$ & $\eta_{\omega}$ & $\tau$\\ 
\hline
Figure~\ref{fig:pipeline} & 0.05 & 0.96 & 0.74 & 88.16 \\
Figure~\ref{fig:pipelineArtifical} & 0.02 & 0.98 & 0.40 & 98.39 \\
Model 1 & 0.05 & 0.96 & 0.34 & 93.38 \\
Model 2 & 0.05 & 0.94 & 0.47 & 89.44 \\
Model 3 & 0.05 & 0.94 & 0.52 & 89.19 \\
Model 4 & 0.20 & 0.98 & 0.76 & 96.93 \\
Model 5 & 0.20 & 0.98 & 0.63 & 96.94 \\
Model 6 & 0.05 & 0.88 & 0.70 & 48.91 \\
Model 7 & 0.02 & 0.98 & 0.78 & 98.78 \\
Model 8 & 0.02 & 0.97 & 0.51 & 98.28 \\
Model 9 & 0.02 & 0.98 & 0.74 & 96.36 \\
Model 10 & 0.02 & 0.97 & 0.37 & 95.03 \\ 
\hline
\end{tabular}
\caption{The quality of obtained meshes:
average quality of elements under $\overline{\eta}$,
the worst element’s quality under $\eta_{\omega}$  and
percentage of elements with quality higher than $0.9$ under $\tau$. Estimation of element's quality adopted from \cite{Remacle:2012}.}	
\label{table_meshQuality}
\end{center}
\end{table}

\section{Acknowledgments}

The present study was carried out in the framework of the research project "Hextreme", funded by the European Research Council (ERC-2015-AdG-694020) and hosted at the Universit\' e catholique de Louvain.

\bibliographystyle{model1-num-names}
\bibliography{sample_arxiv.bib}
\end{document}